\mag=\magstephalf
\input amstex
\documentstyle{amsppt}
\pagewidth{12.6cm}
\pageheight{19.8cm}
\hcorrection{1cm}
\NoBlackBoxes
\def\t#1{\text{\rm #1}}
\def\vep{\varepsilon}
\def\sig{\sigma}
\def\alp{\alpha}
\def\vphi{\varphi}
\def\lam{\lambda}
\def\Lam{\Lambda}
\def\vX{\varXi}
\def\gam{\gamma}
\def\vth{\vartheta}
\def\Del{\Delta}

\def\lot{\overset . \to\otimes}
\def\stcong{\overset . \to\cong}
\def\T{\Tilde}
\def\tr{\operatorname{tr}}

\def\Ch{\operatorname{Ch}}
\def\dim{\operatorname{dim}}
\def\vdim{ \operatorname {\bold{dim}}}
\def\diag{ \operatorname {diag}}
\def\Irr{ \operatorname {Irr}}
\def\id{\t{\rm id}}
\def\E{ \operatorname {End}}
\def\H{ \operatorname {Hom}}

\def\lfl{\lfloor}
\def\rfl{\rfloor}

\def\fS{\frak S}
\def\fq{\frak q}
\def\fgl{\frak g \frak l}

\def\fq{{\frak q}}

\def\pr{{\frak p}}

\def\A{{\Cal A}}
\def\B{{\Cal B}}
\def\C{{\Cal C}}
\def\U{{\Cal U}}
\def\sot{\lot}
\def\pref#1#2#3#4#5#6#7
{\ref\no#1\by#2\paper#3\jour#4\vol#5\pages#6\yr#7\endref}
\topmatter
\title 
A duality of the twisted group algebra of the symmetric group 
and a Lie superalgebra 
\endtitle
\author
Manabu Yamaguchi
\endauthor
\affil
Department of Mathematics \\
Aoyama Gakuin University\\
Chitosedai 6-16-1, Setagaya-ku, Tokoyo, 153 Japan \\
E-mail: yamaguti\@gauss.gem.aoyama.ac.jp\\
\endaffil
\leftheadtext{Manabu Yamaguchi}
\rightheadtext{Duality }
\endtopmatter
\document
\head
0. Introduction
\endhead
The ``character values'' of the irreducible projective 
representations of $\fS_k$, the symmetric group of degree $k$, 
were determined by I\. Schur 
using Schur's $Q$-functions, 
which are indexed by the distinct partitions of $k$, [10], 
in a way analogous to Frobenius' formula 
for the character values of the ordinary irreducible 
representations of $\fS_k$ [2]. 
Behind Frobenius' formula exists a 
duality relation of $\fS_k$ and the general linear group 
$GL(n)$ (the Schur-Weyl duality). 
It is natural to expect the existence of an analogous duality relation between 
the twisted group algebra $\A_k$ (cf\. \thetag{1.2}) 
of $\fS_k$ and some algebra, 
behind Schur's method. 
A\. N\. Sergeev showed that a twisted group algebra $\B_k$ 
(cf\. \thetag{1.3}) of the hyperoctahedral group $H_k$ 
and a Lie superalgebra $\fq(n)$ 
(cf\. \S1, {\bf G}) act on the $k$-th tensor product $W=V^{\otimes k}$ 
of the $2n$-dimensional natural representation $V={\Bbb C}^n\oplus{\Bbb C}^n$ 
of $\fq(n)$, as mutual commutants of each other [11] (in the sense of 
${\Bbb Z}/2{\Bbb Z}$-graded algebras, see \S1, {\bf E}). 
This result motivated our work. 

In this paper, we establish a duality 
relation between $\A_k$ and $\fq(n)$ on a subspace of $W$, and 
give a representation-theoretic explanation of Schur's identity \thetag{1.6} 
adapted to the context of ${\Bbb Z}/2{\Bbb Z}$-graded representations 
by T\. J\'ozefiak (Corollary 4.2). 

In \S3, we construct an isomorphism $\B_k\cong\C_k\lot\A_k$ 
of ${\Bbb Z}/2{\Bbb Z}$-graded algebras (Theorem 3.2), 
where $\C_k$ is the $2^k$-dimensional Clifford algebra 
and $\lot$ denotes the ${\Bbb Z}/2{\Bbb Z}$-graded tensor product 
(cf\. \S1, {\bf E}). 
This isomorphism does imply an embedding $\A_k\hookrightarrow\B_k$, although 
$\A_k$ does not sit in $\B_k$ in an obvious manner (cf. \S1, {\bf D}). 
Then, we give a simple relation between the ${\Bbb Z}/2{\Bbb Z}$-graded 
irreducible representations of $\B_k$ and $\A_k$ (Proposition 3.5). 
Note that J\. R\. Stembridge constructed 
the non-graded simple modules 
of the underlying algebra $|\B_k|$ 
as submodules of non-graded tensor products 
of modules of three twisted group algebras of $H_k$ [13], but 
his description of the simple $|\B_k|$-modules 
does not immediately show a simple relation between $\A_k$ and $\B_k$.

In \S4, we give a submodule $W'$ of $W$, as a simultaneous 
eigenspace of $2^{\lfl k/2\rfl}$ involutions contained in $\C_k$, where 
$\lfl k/2\rfl$ denotes the largest integer not exceeding $k/2$, 
and show that $\A_k$ and $\fq(n)$ act on $W'$ 
as mutual commutants of each other in the context of 
${\Bbb Z}/2{\Bbb Z}$-graded algebras (Theorem 4.1). 

In this paper, all vector spaces, and associative algebras, 
and representations in this paper are assumed to 
be finite dimensional unless otherwise stated, 
and over the complex number field ${\Bbb C}$. 

\subhead
Acknowledgments
\endsubhead
We would like to express our gratitude to 
T\. J\'ozefiak, K\. Koike and I\. Terada 
for various helpful remarks. 

\head
1. Preliminaries and notations
\endhead
\subhead 
A. The hyperoctahedral group $H_k$
\endsubhead
The hyperoctahedral group $H_k$ is the group generated 
by $t$ and the $s_i$, $1\leq i\leq k-1$, subject to relations
$$
\align
&t^2=s_i^2=1\quad(1\leq i\leq k-1),\tag1.1\\
&(s_is_{i+1})^3=1\quad(1\leq i\leq k-2),\quad
(s_is_j)^2=1\quad(|i-j|\geq2),\\
&(ts_i)^2=1\quad(2\leq i\leq k-1),\quad(ts_1)^4=1.
\endalign
$$
$H_k$ is isomorphic to the Weyl group of type $B_k$ or $C_k$, and is also 
sometimes called the group of signed permutations. The subgroup of $H_k$ 
generated by the $s_i$, $1\leq i\leq k-1$, is 
isomorphic to the symmetric group of degree $k$, which we denote by $\fS_k$. 

\subhead
B. Partitions
\endsubhead
Let $P_k$ denote the set of all partitions of $k$ (see [8, p\. 1]), and 
put $P=\coprod_{k\geq0}P_k$. 
For $\lam\in P$, we write $l(\lam)$ for the length of $\lam$, namely 
the number of nonzero parts of $\lam$. Also we write $|\lam|=k$ 
if $\lam\in P_k$.
Let $DP_k$ and $OP_k$ denote the distinct partitions 
(or strict partitions, namely partitions whose parts are distinct) 
and the odd partitions (namely partitions whose parts are all odd) 
of $k$ respectively. 
Let $DP^{+}_k$ and $DP^{-}_k$ be the sets of all $\lam\in DP_k$ 
such that $(-1)^{k-l(\lam)}=+1$ and $-1$ respectively. Note that 
$(-1)^{k-l(\lam)}$ equals the signature of permutations with 
cycle type $\lam$. 
We also put $DP=\coprod_{k\geq0}DP_k$ and $OP=\coprod_{k\geq0}OP_k$. 
Note that these notations, $DP$ and $OP$, were used by Stembridge 
in [12], [13]. 

\subhead
C. The subring $\Omega$ of the ring of the symmetric functions
\endsubhead
Let $\Lam$ denote the ring of the symmetric functions with 
coefficients in ${\Bbb C}$. 
Note that our $\Lam$ is the scalar extension of the $\Lam$ in [8], which 
is a ${\Bbb Z}$-algebra, to ${\Bbb C}$. 
We have $\Lam=\bigoplus_{k\geq0}\Lam^k$, 
where $\Lam^k$ denotes the subspace of $\Lam$ consisting of the homogeneous 
elements of degree $k$. 

For each $r\geq1$, let $p_r$ denote the $r$-th power sum 
$\sum_{i\geq1}x_i^r\in\Lam$. 
The $p_r$, $r\geq 1$, are algebraically independent over ${\Bbb C}$ and 
we have $\Lam={\Bbb C}[p_1,p_2,\dots]$. 
For each partition $\lam=(\lam_1,\lam_2,\dots,\lam_l)\in P$, 
let $p_\lam=p_{\lam_1}p_{\lam_2}\cdots p_{\lam_l}$. Then 
$\{p_{\lam}\,;\,\lam\in P\}$ is a basis of $\Lam$ 
(cf\. [8, \S{I}, sect\. 2.]). 

Let $\Omega$ denote the subring of $\Lambda$ 
generated by the power sums of odd degrees, namely the $p_r$, 
$r=1,3,5$, $\dots$ . 
Put $\Omega^k=\Omega\cap\Lam^k$. 
Then we have $\Omega=\bigoplus_{k\geq0}\Omega^k$, 
and $\{p_{\lam}\,|\,\lam\in OP\}$ is a basis of $\Omega$. 

For $\lam\in DP$, let $Q_{\lam}\in \Lam$ 
denote Schur's $Q$-function indexed by $\lam$ 
(cf\. [10], [12, sect\. 6]). 
Then $\{Q_{\lam}\,|\,\lam\in DP\}$ is also a basis of $\Omega$. 

\subhead
D. Twisted group algebras $\A_k$ and $\B_k$
\endsubhead
For each $k\geq1$, let $\A_k$ denote the associative algebra generated by 
the elements $\gamma_i$, $1\leq i\leq k-1$, subject to relations
$$
\align
&\gamma_i^2=-1\;\;(1\leq i\leq k-1),\;\;
(\gamma_i\gamma_{i+1})^3=-1\;\;(1\leq i\leq k-2),\tag1.2\\
&(\gamma_i\gamma_j)^2=-1\;\;(|i-j|\geq2).
\endalign
$$ 
$\A_k$ is isomorphic to a twisted group algebra of $\fS_k$ 
(see below). 
We regard $\A_k$ as a ${\Bbb Z}_2$-graded algebra by giving 
degree $1\in{\Bbb Z}_2$ to the generators $\gamma_i$ 
for all $1\leq i\leq k-1$. 
In the following, we abbreviate ${\Bbb Z}/2{\Bbb Z}$ as ${\Bbb Z}_2$. 

For each $k\geq1$, let $\B_k$ denote the associative algebra generated by 
$\tau$ and the $\sig_i$, $1\leq i\leq k-1$, subject to relations 
$$
\align
&\tau^2=\sig_i^2=1\quad(1\leq i\leq k-1),
\quad(\sig_i\sig_{i+1})^3=1\quad(1\leq i\leq k-2),
\tag1.3\\
&(\sig_i\sig_j)^2=1\quad(|i-j|\geq2),\quad
(\tau\sig_i)^2=1\quad(2\leq i\leq k-1),\\
&(\tau\sig_1)^4=-1.\\
\endalign
$$ 
$\B_k$ is isomorphic to a twisted group algebra of $H_k$ 
(again see below). 
We regard $\B_k$ as a ${\Bbb Z}_2$-graded algebra by giving 
degree $1$ to the generator $\tau$ and degree $0$ to 
the generator $\sig_i$ for all $1\leq i\leq k-1$. 

\remark{Remark} The structures of $H^2(\fS_k,{\Bbb C}^{\times})$ 
and $H^2(H_k,{\Bbb C}^{\times})$ were clarified by I\. Schur [10] and 
J\. W\. Davies, A\. O\.  Morris [1] respectively, as follows: 
$$
H^2(\fS_k,{\Bbb C}^{\times})\cong
\left\{\matrix
\format\l&\l\\
0         &\hskip5mm\t{\rm for\;\;}k\leq3\\
{\Bbb Z}/2{\Bbb Z}&\hskip5mm\t{\rm for\;\;}k\geq4\;\;
\endmatrix\right.\tag1.4 
$$

$$
H^2(H_k, {\Bbb C}^{\times})\cong\left\{
\matrix
\format\l&\quad\l\\
0&\t{if $k=1$}\\
{\Bbb Z}/2{\Bbb Z}&\t{if $k=2$}\\
{\Bbb Z}/2{\Bbb Z}\times{\Bbb Z}/2{\Bbb Z}
&\t{if $k=3$}\\
{\Bbb Z}/2{\Bbb Z}\times{\Bbb Z}/2{\Bbb Z}\times{\Bbb Z}/2{\Bbb Z}
&\t{if $k\geq4$}.\\
\endmatrix\right.\tag1.5
$$
If $k\geq4$, then any twisted group algebra ${\Bbb C}^{\alp}\fS_k$ 
of $\fS_k$ with a non-trivial $2$-cocycle $\alp$ is isomorphic to $\A_k$ 
(cf\. [10], [12, Lem\. 1.1]). If $k\leq3$, then $\A_k$ is isomorphic to 
the ordinary group algebra ${\Bbb C}\fS_k$, via $\gam_i\mapsto\sqrt{-1}s_i$. 

Moreover, if $k\geq2$, then $\B_k$ is isomorphic to 
a twisted group algebra ${\Bbb C}^{\alp}H_k$ of $H_k$ with a 
non-trivial $2$-cocycle $\alp$ (cf\. [13, Prop\. 1.1]). 
Note that the cocycle associated with $\B_k$ does not restrict 
to the cocycle associated with $\A_k$. 
\endremark

In the following, we consistently use the formulation of the 
${\Bbb Z}_2$-graded representations, namely 
the ${\Bbb Z}_2$-graded modules of ${\Bbb Z}_2$-graded algebras 
(superalgebras), as was used in [3] and [4]. 
This formulation is slightly different 
from the traditional parametrization of the non-graded 
representations, but there exists a explicit relation between 
the ${\Bbb Z}_2$-graded representations 
of a ${\Bbb Z}_2$-graded algebra $A$ 
and the non-graded representations of the underlying algebra $|A|$ 
(cf\. [3, Lem\. (2.8), Cor\. (2.16), Prop\. (2.17)], 
[4, Prop\.  2.5, Cor\.  2.6]). 

\subhead
E. Semisimple superalgebras
\endsubhead
This theory was developed by T\. J\'ozefiak in [3], 
which we mostly follow. 
A ${\Bbb Z}_2$-graded algebra $A$, which is called 
a {\bf superalgebra} in this paper, 
is called {\bf simple} if it has no ${\Bbb Z}_2$-graded two-sided 
ideals except itself and $0$. 

Let $V$ be a ${\Bbb Z}_2$-graded vector space, namely a vector space 
with a fixed direct sum decomposition $V=V_0\oplus V_1$. We write 
$\vdim V$ for the pair $(\dim V_0,\dim V_1)$. If 
$W=W_0\oplus W_1$ is another ${\Bbb Z}_2$-graded vector space, then 
the vector space $\H(V,W)$ consisting of all linear maps from 
$V$ to $W$ has a ${\Bbb Z}_2$-gradation defined as follows:
$$
\H^{\alp}(V,W)=\{f\in\H(V,W)\,;\,f(V_{\beta})\subset W_{\alp+\beta} 
\t{ for all $\beta\in {\Bbb Z}_2$}\}
$$
for each $\alp\in {\Bbb Z}_2$. In particular, it can be 
easily checked that the endomorphism algebra $\E(V)=\H(V,V)$, 
which is isomorphic to the full matrix ring $M_{n+m}$ where $\vdim V=(n,m)$, 
can be regarded as a superalgebra with this gradation. 
This superalgebra is denoted by $M(n,m)$. 
It is clear that $M(n,m)$ is simple, since the underlying algebra 
$|M(n,m)|=M_{n+m}$ is a simple algebra. 

There exist another type of simple superalgebras. 
Let $Q(n)$ denote a subsuperalgebra of $M(n,n)$ 
which consists of all $2n\times 2n$-matrices of the form 
$\pmatrix
C&D\\
D&C
\endpmatrix$, 
where $C$ and $D$ are $n\times n$-matrices. 
It is easy to show that $Q(n)$ is a simple superalgebra. 

\proclaim{Theorem 1.1}\t{\rm(cf\. [3, Th\. 2.6], [4, Th\. 2.1], [15])} 
A simple superalgebra is isomorphic to either 
$M(n,m)$ for some $n,m$ or $Q(n)$ for some $n$.
\endproclaim

A simple superalgebra $A$ is said to be 
of {\bf type} $M$ (resp\. of {\bf type} $Q$) if 
$A$ is isomorphic to $M(n,m)$ for some $n,m$ 
(resp\. is isomorphic to $Q(n)$ for some $n$). 

Let $V$ be an $A$-{\bf module}, namely 
a ${\Bbb Z}_2$-graded vector space 
$V=V_0\oplus V_1$ together with an algebra homomorphism 
$\rho\:A\to\E(V)$, which preserves ${\Bbb Z}_2$-gradations. 
Then we call $\rho$ a {\bf representation} of $A$ in $V$, and 
simply write $\rho(a)v=av$ for all $a\in A$ and $v\in V$. 
A ${\Bbb Z}_2$-graded subspace $W$ of $V$ is called an 
$A$-{\bf submodule} if it is stable  under $\rho(A)$. 
We say that $V$ is {\bf simple} if 
it has no $A$-submodules except itself and $0$. 

Let $V$ and $W$ be two $A$-modules. 
For each $\alp\in{\Bbb Z}_2$, 
let $\H^{\alp}_A(V,W)$ denote 
the subspace of $\H^{\alp}(V,W)$ consisting 
of all elements $f\in \H^{\alp}(V,W)$ such that 
$f(av)=(-1)^{\alp\cdot\beta}af(v)$ 
for all $a\in A_{\beta}$ $(\beta\in{\Bbb Z}_2)$, $v\in V$. 
Put $\H^{\cdot}_A(V,W)=\H^{0}_A(V,W)\oplus\H^{1}_A(V,W)$ and put 
$\E^{\cdot}_A(V)=\H^{\cdot}_A(V,V)$. It is clear that 
$\E^{\cdot}_A(V)$ is a subsuperalgebra of $\E(V)$ with this gradation. 
We call $\E^{\cdot}_A(V)$ the {\bf supercentralizer } of 
$\rho(A)$ or $A$ in $\E(V)$, where $\rho$ is the representation 
of $A$ associated with $V$. 
The {\bf shift} of $V$, denoted by $\overline{V}$ in this paper, is defined 
to be the same vector space as $V$ with the switched grading, namely 
$\overline{V}_0=V_1$ and $\overline{V}_1=V_0$, and 
with the homomorphism $\overline{\rho}\:A\to\E(\overline{V})$ 
defined by $\overline{\rho}(a)=(-1)^{\alp}\rho(a)$ for $a\in A_{\alp}$, 
where $\rho$ is the superalgebra homomorphism associated with $V$. 
 Note that we have 
$\H^{0}_A(\overline{V},W)=\H^{1}_A(V,W)$ and 
$\H^{1}_A(\overline{V},W)=\H^{0}_A(V,W)$. 
Two $A$-modules $V$ and $W$ are called {\bf isomorphic} (resp\. 
{\bf strictly isomorphic}) if there exists an invertible linear map 
$f\in\H^{\cdot}_A(V,W)$ 
(resp\. $f\in\H^{0}_A(V,W)$). 
If this is the case, we write $V\cong_AW$ 
(resp\. $V\stcong_AW$). In view of the following theorem, 
if $V$ and $W$ are simple $A$-modules, 
we have $V\cong_AW$ if and only if $V\stcong_AW$ or $\overline{V}\stcong_AW$.

The following theorem contains an analogue of Schur's Lemma.
\proclaim{Theorem 1.2}\t{\rm(cf\. [3, Prop\. 2.17],
[4, Prop\.  2.5, Cor\.  2.6])} 
Let $V=V_0\oplus V_1$ be a simple $A$-module. 

\t{\rm(1)} The supercentralizer $\E^{\cdot}_A(V)$ is isomorphic to 
$M(1,0)\cong{\Bbb C}$, if and only if 
$V_0$ is not isomorphic to $V_1$ as an $A_0$-module. 
\vskip0.3cm

\t{\rm(2)} The supercentralizer $\E^{\cdot}_A(V)$ 
is isomorphic to $Q(1)\cong\C_1$, 
if and only if $V_0$ is isomorphic to $V_1$ as an $A_0$-module. 
\endproclaim
If a simple $A$-module $V$ satisfies either of the equivalent 
conditions in (1) (resp\. (2)), 
then we say that it is 
of {\bf type} $M$ (resp\. of {\bf type} $Q$). 
If $V$ is a simple $A$-module of type $M$, then so is $\overline{V}$ and 
we have $V \not \stcong_A \overline{V}$.  
If $V$ is a simple $A$-module of type $Q$, then so is $\overline{V}$ and 
we have $V \stcong_A \overline{V}$. In both cases we have 
$V\cong_A\overline{V}$. 

Let $A$ and $B$ be two superalgebras and let 
$V$ (resp\. $W$) be an $A$ (resp\. $B$)-module. 
The vector space $A\otimes B$ can be turned into a superalgebra, 
where the grading is defined by $(A\otimes B)_{\alp}=
\oplus_{\beta+\gamma=\alp}A_{\beta}\otimes B_{\gamma}$ for each 
$\alp\in{\Bbb Z}_2$, and the multiplication is defined by 
$$
(a\otimes b)(c\otimes d)=(-1)^{\beta\cdot\gam}ac\otimes bd
$$
for any $a\in A$, $b\in B_{\beta}$, $c\in A_{\gam}$, and $d\in B$ 
($\beta$, $\gam\in{\Bbb Z}_2$). 
This superalgebra is called the {\bf supertensor product} 
of $A$ and $B$ and denoted by $A\lot B$. 
Note that we have an isomorphism of superalgebras 
$\omega_{A,B}\:A\lot B\to B\lot A$ determined by 
$\omega_{A,B}(a\lot b)=(-1)^{\alp\cdot\beta}b\lot a$ for 
all homogeneous elements $a\in A_{\alp}$ and $b\in B_{\beta}$ 
($\alp$, $\beta\in{\Bbb Z}_2$). 
The following states that the supertensor product of simple superalgebras 
is again a simple superalgebra. 
\proclaim{Theorem 1.3} \t{\rm [3, Prop\. 2.10], [15]} 
There exist isomorphisms of superalgebras 
\roster
\item"(a)" $M(r,s)\lot M(p,q)\cong M(rp+sq,rq+sp)$,
\item"(b)" $M(r,s)\lot Q(n)\cong Q(rn+sn)$,
\item"(c)" $Q(m)\lot Q(n)\cong M(mn,mn)$.
\endroster
\endproclaim

Moreover, $V\otimes W$ can be regarded as a $A\lot B$-module as follows:
$$
(a\lot b)(v\otimes w)=(-1)^{\beta\cdot\alp}av\otimes bw
$$
for any $a\in A$, $b\in B_{\beta}$, $v\in V_{\alp}$, and $w\in W$ 
($\alp$, $\beta\in{\Bbb Z}_2$), and 
$(V\otimes W)_{\alp}
=\oplus_{\beta+\gamma=\alp}V_{\beta}\otimes W_{\gamma}$ 
for $\alp\in{\Bbb Z}_2$. 
This $A\lot B$-module is called the 
{\bf supertensor product} of $V$ and $W$ and denoted by $V\lot W$. 
Note that the supertensor product is symmetric, namely 
the $B\lot A$-module obtained from the $A\lot B$-module 
$V\lot W$ via $\omega_{B,A}$ is isomorphic to $W\lot V$ by the map 
determined by $v\lot w\mapsto (-1)^{\alp\cdot\beta}w\lot v$ for 
all homogeneous $v\in V_{\alp}$ and $w\in W_{\beta}$ 
($\alp$, $\beta\in{\Bbb Z}_2$). 

The following follows from Schur's lemma.
\proclaim{Theorem 1.4} 
 Let $A$, $B$ be superalgebras, and put $C=A\lot B$. 
 Let $U$ (resp\. $W$) be a simple $A$-module (resp\. simple $B$-module), and 
put $V=U\lot W$. 
\roster 
\item"(a)" If $U$, $W$ are of type $M$, then $V$ is a simple $C$-module 
of type $M$. We have $\overline U \lot W \stcong_C U \lot \overline W 
\stcong_C \overline V$, $\overline U \lot \overline W \stcong_C V$ 
and $V|_{A}\stcong_A U^{\oplus k'} \oplus \overline U^{\oplus l'}$, 
$V|_{B}\stcong_B W^{\oplus k} \oplus \overline W^{\oplus l}$ 
where $(k,l)=\vdim U$, $(k',l')=\vdim W$. 
\item"(b)" If one of $U$ and $W$ is of type $M$ and the other 
is of type $Q$, then from the symmetry of the supertensor product 
we have only to state for the case where 
$U$ is of type $M$ and $W$ is of type $Q$. 
Then $V$ is a simple $C$-module of type $Q$. 
We have $\overline U \lot W \stcong_C V$ and 
$V|_{A}\stcong_A U^{\oplus n'} \oplus \overline U^{\oplus n'}$, 
$V|_{B}\stcong_B W^{\oplus k+l}$ where 
$(k,l)=\vdim U$, $(n',n')=\vdim W$. 
\comment
\item"(c)" If $U$ is of type $Q$, $W$ is of type $M$, 
then $V$ is a simple $C$-module of type $Q$. 
We have $U \lot \overline W \stcong_C V$ 
and $V|_{A}\stcong_A U^{\oplus k'+l'}$ 
where $(k',l')=\vdim W$. 
\endcomment
\item"(c)" Suppose $U$ and $W$ are of type $Q$. 
Fix $x\in \E^{1}_A(U)$ (resp\. $y\in \E^{1}_B(W)$) satisfying 
$x^2=-1$ (resp\. $y^2=-1$). Then $x\lot y\in \E^{0}_C(U\lot W)$ satisfies 
$(x\lot y)^2=-1$. Let $V^{\pm}$ be the $(\pm\sqrt{-1})$-eigenspace 
of $x\lot y$ respectively. 
Then $V=V^{+}\oplus V^{-}$ is a decomposition into simple $C$-modules, 
both of type $M$. 
We have $\overline{V^{+}} \stcong_C V^{-}$ 
and $V^{+}|_{A}\stcong_A V^{-}|_{A}\stcong_A U^{\oplus n'}$, 
$V^{+}|_{B}\stcong_B V^{-}|_{B}\stcong_B U^{\oplus n}$ 
where $(n,n)=\vdim U$, $(n',n')=\vdim W$. 
\endroster

Moreover, the above construction gives all simple $A\lot B$-modules. 
\endproclaim

Let $\Irr A$ denotes the set of all isomorphism classes 
(not strict isomorphism classes) of simple $A$-modules for any superalgebra 
$A$.
\proclaim{Corollary 1.5} We have a bijection $\Irr A\lot B 
\overset \sim \to \to \Irr A \times \Irr B$. 
\endproclaim

For any $A$-module $V$, the following conditions are equivalent. 
\roster
\item $V$ is a sum of simple $A$-submodules, 
\item $V$ is a direct sum of simple $A$-submodules, 
\item For any $A$-submodule $W$ of $V$, there exists 
a $A$-submodule $W'$ of $V$ such that $V=W \oplus W'$. 
\endroster 
We call an $A$-module $V$ {\bf semisimple} if it satisfies one of 
the above equivalent conditions (1)-(3). 
 
For any superalgebra $A$, the following conditions are equivalent 
(cf\. [3, Prop\. 2.4], [3, Cor\. 2.12], [4, Th\. 2.2, Cor\. 2.3]): 
\roster
\item $A$ is a semisimple (regular) $A$-module, namely 
$A$ is a direct sum of simple $A$-submodules, 
\item every $A$-module is semisimple, 
\item $A$ is a direct sum of simple superalgebras. 
\endroster 
We call a superalgebra $A$ {\bf semisimple} if it satisfies one of 
the above equivalent conditions (1)-(3). 
Note that the supertensor product of two semisimple superalgebras 
is semisimple. 
\proclaim{Theorem 1.6}\t{\rm(cf\. [3, Cor\. 2.12], 
[4, Th\. 2.2, Cor\. 2.3])} Let $A$ be a semisimple superalgebra. 

\t{(1)} There exist integers $m,q\geq0$ and $k_i,l_i\geq0$ 
$(1\leq i\leq m)$ and $n_j\geq 0$ $(1\leq j\leq q)$ such that 
$$
A\cong 
\bigoplus_{i=1}^{m}M(k_i,l_i)\oplus
\bigoplus_{j=1}^{q}Q(n_j)
$$
as a superalgebra. The following data are determined 
by $A$, and conversely the following data determine $A$: 
{\rm(a)} $m=m(A)$, {\rm(b)} $q=q(A)$, {\rm(c)} the multiset of the unordered 
pairs $\{k_i,l_i\}$, and {\rm(d)} the multiset of the numbers $n_j$. 

\t{\rm(2)} The number of the isomorphism 
classes of the simple $A$-modules is equal to $m(A)+q(A)$.
\endproclaim

Note that the set of simple direct summands of $A$, isomorphic to 
$M(k_i,l_i)$ $(1\leq i\leq m)$ and $Q(n_j)$ $(1\leq j\leq q)$, 
is uniquely determined by $A$. 
We call them the {\bf simple components} of $A$. 
We should also remark the equivalence of the semisimplicity 
of a superalgebra to that of the underlying algebra:
\proclaim{Proposition 1.7} \t{(cf\. [3, Cor\. 2.16])} A superalgebra 
$A$ is semisimple if and only if an underlying algebra $|A|$ is semisimple. 
\endproclaim

Let $A$ be a semisimple superalgebra, and let $V$ be an $A$-module. 
Then there exist a finite number of simple $A$-submodules 
$U_1,U_2,\dots,U_l$ of $V$ such that 
$V=U_1\oplus U_2\oplus\cdots\oplus U_l$ as an $A$-module. 
For each $1\leq r\leq l$, put 
$U(V)_r=\sum_{U_s\cong_AU_r}V_s$. Then we can choose a subset $R$ 
of $\{1,2,\dots,l\}$ such that $V=\bigoplus_{r\in R}U(V)_r$.  
Note that the set of $A$-submodules $U(V)_r$, $r\in R$, 
is uniquely determined by $V$. The $U(V)_r$ are called the 
$A$-{\bf homogeneous components} of $V$. 

In the next section, \S2, we will precisely state 
the double centralizer theorem for 
semisimple superalgebras, in which we will show that there exists 
a $1-1$ correspondence between the simple $A$-modules and the simple 
$\E^{\cdot}_A(V)$-modules appearing in an $A$-module $V$. 

\subhead
F. Character formulas for $\A_k$ and $\B_k$
\endsubhead
Let $A$ be a superalgebra: $A=A_0\oplus A_1$, 
and let $V$ be an $A$-module. 
Define a linear map $\Ch[V]\:A\to{\Bbb C}$ by 
$\Ch[V](a)=\tr_V(a)$ for all $a\in A$. It is called the {\bf character} 
of the $A$-module $V$. $\Ch[V]$ only depends on the isomorphism class of $V$. 
We also say that $V$ {\bf affords} the character $\Ch[V]$. 
Let $\{U_r\}_{r\in R}$ be a complete set of representatives 
of the isomorphism classes of the simple $A$-modules. 
Let $Z_0(A^{*})$ denote the subspace of 
$A^{*}=\H(A,{\Bbb C})$ consisting of all elements 
$f\:A\to{\Bbb C}$ such that $f(ab)=f(ba)$ for all $a$, $b\in A$ and 
such that $f|_{A_1}=0$. The character of any $A$-module belongs to 
$Z_0(A^{*})$. If $A$ is semisimple, then 
$\{\Ch[U_r]\,;\, r\in R\}$ is a basis of $Z_0(A^{*})$. In this case, 
the isomorphism classes of $A$-modules are uniquely determined by their 
characters.

Let $\fS'_k$ be the subgroup of $(\A_k)^{\times}$ 
generated by $-1,\gam_1,\dots,\gam_{k-1}$ ($-1$ is not necessary if 
$k\geq 4$). 
Then $\fS'_k$ is a double cover (a central extension with 
a ${\Bbb Z}_2$ kernel) of $\fS_k$ with a 
group homomorphism $\pi\:\fS'_k\to\fS_k$ 
defined by $\pi(-1)=1$ and $\pi(\gam_j)=s_j$ for all $1\leq j\leq k-1$. 
If $k\geq 4$, then $\fS'_k$ is a representation group of $\fS_k$. 
Let $\gam^{\mu}$ be an element of 
$\fS'_k$ defined by taking a particular reduced expression of 
a permutation of cycle type $\mu$ and replacing each $s_i$ with $\gam_i$, 
namely: 
$$
\gam^{\mu}=
(\gam_1\cdots\gamma_{\mu_1-1})
(\gam_{\mu_1+1}\cdots\gam_{\mu_1+\mu_2-1})\cdots
(\gam_{\mu_1+\cdots+\mu_{r-1}+1}\cdots
\gam_{\mu_1+\cdots+\mu_r-1})
$$
where $r=l(\mu)$. 
Let $\vphi$ be a character of an $\A_k$-module. 
Then we have $\vphi(\gam^{\mu})=0$ unless $\mu\in OP_k$ 
\comment
The characters of $\A_k$ are completely determined 
by their values at the elements $\gam^{\mu}$, $\mu\in OP_k$,
\endcomment
(cf\. [10], [12]). 
I\. Schur explicitly described the characters 
of the simple $\A_k$-modules. 
We review Schur's result in the form translated by 
T\. J\'ozefiak into the language of ${\Bbb Z}_2$-graded representations [4]. 
Let $\vep\:DP_k\to{\Bbb Z}_2$ be defined by 
$\vep(\nu)=0$ if $\nu\in DP^{+}_k$ and $\vep(\nu)=1$ if 
$\nu\in DP^{-}_k$. 

\proclaim{Theorem 1.8} \t{[4], [10]} Define 
the $\vphi_{\nu}\in Z_0((\A_k)^{*})$, $\nu\in DP_k$, by 
$$
(\sqrt{2})^{l(\mu)}p_{\mu}=\sum_{\nu\in DP_k}
\vphi_{\nu}(\gam^{\mu})
(\sqrt{2})^{-l(\nu)-\vep(\nu)}Q_{\nu}\tag1.6
$$
for all $\mu\in OP_k$, 
where $Q_{\nu}$ denotes the Schur $Q$-function 
(cf\. \t{[10], [12]}).
Then the $\vphi_{\nu}$, $\nu\in DP_k$, give all 
irreducible characters of $\A_k$. 

For each $\nu\in DP_k$, fix a simple $\A_k$-module $V_{\nu}$ 
which affords the irreducible character $\vphi_{\nu}$ of $\A_k$.
Then $V_{\nu}$ is of type $M$ (resp\. type $Q$) 
if $\nu\in DP_k^{+}$ (resp\. $\nu\in DP_k^{-}$). 
\endproclaim

Let $H'_k$ be the subgroup of $(\B_k)^{\times}$ 
generated by $-1$, $\tau$, $\sig_1,\dots,\sig_{k-1}$ ($-1$ is not necessary 
if $k\geq 2$). 
Then $H'_k$ is a double cover of $H_k$ with a group 
homomorphism $\pi\:H'_k\to H_k$ 
defined by $\pi(-1)=1$, $\pi(\tau)=t$, 
and $\pi(\sig_j)=s_j$ for all $1\leq j\leq k-1$. 
Let $\sig^{(\lam,\mu)}$ denote a similarly defined element of $H'_k$ mapping to a representative of the conjugacy class of $H_k$ indexed by $(\lam,\mu)$: 
$$
\align
&\sig^{(\lam,\mu)}=x_1x_2\cdots y_1y_2\cdots,\\
&x_i=\sig_{\lam_1+\cdots+\lam_{i-1}+1}\cdots\sig_{\lam_1+\cdots+\lam_i-1},\\
&y_i=\sig_{|\lam|+\mu_1+\cdots+\mu_{i-1}+1}\cdots
\sig_{|\lam|+\mu_1+\cdots+\mu_i-1}\tau_{|\lam|+\mu_1+\cdots+\mu_i},\\
&\tau_j=\sig_{j-1}\cdots\sig_1\tau\sig_1\cdots\sig_{j-1}\;\;.
\endalign
$$
Let $\psi$ be a character of a $\B_k$-module. 
Then we have $\psi(\sig^{(\lam,\mu)})=0$ unless $\lam\in OP_k$ and $\mu=\phi$ 
\comment
The characters of $\B_k$ are completely determined 
by their values at the elements $\sig^{(\mu,\phi)}$, $\mu\in OP_k$, 
\endcomment
(cf\. [1], [13]). 
A\. N\. Sergeev showed the following explicit formula for the characters 
of the simple $\B_k$-modules using his duality relation 
which we will review in the next subsection, \S1, {\bf G}. 
Define a map $d\:P_k\to{\Bbb Z}_2$ by 
$d(\nu)=1$ if $l(\nu)$ is odd, and $d(\nu)=0$ if $l(\nu)$ is even.

\proclaim{Theorem 1.9} Define the functions  
$\psi_{\nu}\in Z_0((\B_k)^{*})$, $\nu\in DP_k$, by 
$$
2^{l(\mu)}p_{\mu}=
\sum_{\nu\in DP_k}\psi_{\nu}(\sig^{(\mu,\phi)})
(\sqrt{2})^{-l(\nu)-d(\nu)}Q_{\nu}
\tag1.7
$$
for all $\mu\in OP_k$. 
Then the $\psi_{\nu}$, $\nu\in DP_k$, give all 
irreducible characters of $\B_k$. 

For each $\nu\in DP_k$, fix a simple $\B_k$-module 
$W_{\nu}$ which affords the character $\psi_{\nu}$. 
Then $W_{\nu}$ is of type $M$ (resp\. type $Q$) 
if $l(\nu)$ is even (resp\. odd).
\endproclaim

\subhead
G. Sergeev's duality relation 
\endsubhead
First we introduce the Lie superalgebra $\fq(n)$, 
sometimes called the queer Lie superalgebra. 
A standard reference for Lie superalgebras is [7]. 

Let $\fgl(n/m)$ (denoted by $l(n,m)$ in [7]) denote the Lie superalgebra 
whose underlying vector space is that of (the superalgebra) 
$M(n,m)$ and the Jacobi product 
$[\;\;,\;\;]\:\fgl(n/m)\times \fgl(n/m)\to\fgl(n/m)$ is defined by 
$[X,Y]=XY-(-1)^{\alp\cdot\beta}YX$ 
for all $X\in \fgl(n/m)_{\alp}$ and $Y\in \fgl(n/m)_{\beta}$ 
$(\alp, \beta\in{\Bbb Z}_2)$. 
Let $\fq(n)$ denote the Lie superalgebra for the superalgebra 
$Q(n)$ in the same manner. 
Let $\U_n=\U(\fq(n))$ denote the universal enveloping algebra of 
$\fq(n)$, which can be regarded as a superalgebra 
(cf\. [7]). 

The Lie superalgebra $\fq(n)$ naturally acts on a $2n$-dimensional 
space $V$ with a fixed ${\Bbb Z}_2$-gradation 
$V=V_0\oplus V_1$, $\dim{V_0}=\dim{V_1}=n$. 
Therefore $\U_n$ also acts on $V$. 
Let $W$ denote the $k$-fold supertensor product $W=V^{\lot k}$ of $V$. 
We define a representation $\Theta\:\U_n\to\E(W)$ by 
$$
\Theta(X)
(v_1\lot\cdots\lot v_k)
=\sum_{j=1}^{k}(-1)^{\alp\cdot(\beta_1+\cdots+\beta_{j-1})}
v_1\lot\cdots\lot\overset j \to{\overset\smile\to{Xv_j}}\lot\cdots\lot v_k
\tag1.8
$$
for all $X\in\fq(n)_{\alp}$ and $v_i\in V_{\beta_i}$ $(1\leq i\leq k)$ 
and $\alp$, $\beta_i\in{\Bbb Z}_2$ $(1\leq i\leq k)$. 
Actually this action can be defined using a superalgebra homomorphism 
$\Del\:\U_n\to\U_n\lot\U_n$ called the coproduct, determined by 
$\Del(X)=1\lot X+ X\lot1$ for all $X\in \fq(n)$. 
Namely, if we put 
$\Del^{(k)}=(\oversetbrace{k-2}\to{\id \lot\cdots\lot \id }\lot \Del)
\circ \cdots \circ (\id \lot \Del) \circ \Del\:\U_n\to 
\U_n^{\lot k}$, then we have 
$\Theta(X)(v_1\lot\cdots\lot v_k)
=\Del^{(k)}(X)(v_1\lot\cdots\lot v_k)$. 
Note that $\U_n$ is an infinite dimensional superalgebra. 
However, for a fixed number $k$, 
$\U_n$ acts on $W$ through its finite dimensional image in $\E(W)$. 
Therefore we can use the results in \S1, {\bf E} 
on finite dimensional superalgebras 
and their finite dimensional modules. 

Next define a representation $\Psi\:\B_k\to\E(W)$ by
$$
\align
&\Psi(\tau)(v_1\lot\cdots\lot v_k)
=(Pv_1)\lot v_2\lot\cdots\lot v_k\tag1.9\\
&\Psi(\sig_i)(v_1\lot\cdots\lot v_k)
=(-1)^{\beta_i\cdot\beta_{i+1}}
v_1\lot\cdots\lot v_{i+1}\lot v_i\lot\cdots\lot v_k\\
&\hskip4cm
(v_i\in V_{\beta_i},\;\beta_i\in {\Bbb Z}_2,\;1\leq i\leq k-1)
\endalign
$$
where $P=\pmatrix
0&-\sqrt{-1}I_n\\
\sqrt{-1}I_n&0
\endpmatrix\in M(n,n)_1$. 

Let $W'$ be a $\U_n$-submodule of $W$. 
Since $\fq(n)_0\cong\fgl(n,{\Bbb C})$ 
as a Lie algebra, and $V$ is a sum 
of two copies of natural representations of $\fgl(n,{\Bbb C})$ 
($V_0$ and $V_1$), 
this embeds $W'|_{\fq(n)_0}$ into a sum of tensor powers of the natural 
representation, so that this representation of $\fgl(n,{\Bbb C})$ can 
be integrated to a polynomial representation $\theta_{W'}$ of 
$GL(n,{\Bbb C})$. Let $\Ch[W']$ denote the character of $\theta_{W'}$, 
namely $\Ch[W'](x_1,x_2,\dots,x_n)$ is the trace of 
$\theta_{W'}(\diag(x_1,x_2,\dots,x_n))$. 

We review Sergeev's result. 
He established a duality relation between $\B_k$ and $\U_n$ using 
the double centralizer theorem for semisimple superalgebras, which we will 
precisely state in the next section, \S2. 

\proclaim{Theorem 1.10} \t{\rm[11]} 
\t{\rm(1)} Two superalgebras 
$\B_k$ and $\U_n$ act on $W$ as mutual supercentralizers of each other: 
$$
\E^{\cdot}_{\Theta(\U_n)}(W)=\Psi(\B_k),
\quad \E^{\cdot}_{\Psi(\B_k)}(W)=\Theta(\U_n). 
$$

\t{\rm(2)} The simple $\B_k$-module $W_{\nu}$ $(\nu\in DP_k)$ 
occurs in $W$ if and only if $l(\nu)\leq n$.

\t{\rm(3)} For each $\nu\in DP_k$, let $U_{\nu}$ denote 
a simple $\U_n$-module corresponding to $W_{\nu}$ in $W$, 
in the sense of (the final part of) Theorem 2.1. 
Then it follows that 
$$
\Ch[U_{\nu}](x_1,x_2,\dots,x_n)
=(\sqrt{2})^{d(\nu)-l(\nu)}Q_{\nu}(x_1,x_2,\dots,x_n).\tag1.10
$$
\endproclaim 

For each $\nu\in DP_k$ such that $l(\nu)$ is odd, 
fix $y_{\nu}\in\E^{1}_{\B_k}(W_{\nu})$ 
such that $y_{\nu}^2=-1$ 
and let $u_{\nu}\in\E^{1}_{\U_n}(U_{\nu})$ 
be defined using $y_{\nu}$ as in \thetag{2.2} 
in which 
$y_i$ and $x_i$ are replaced by $u_{\nu}$ and $y_{\nu}$ respectively. 
Then we have 
$$
W\stcong_{\B_k\lot \U_n}\bigoplus_{l(\nu):\t{\rm\;even}}
W_{\nu}\lot U_{\nu}\oplus
\bigoplus_{l(\nu):\t{\rm\;odd}}
(W_{\nu}\lot U_{\nu})^{+}\tag1.11
$$
where $(W_{\nu}\lot U_{\nu})^{\pm}$ denotes the simple 
$\B_k\lot\U_n$-module obtained as the $\pm\sqrt{-1}$-eigenspace of 
$y_{\nu}\lot u_{\nu}\in\E^{0}_{\B_k\lot\U_n}(W_{\nu}\lot U_{\nu})$ 
respectively. 

\head
2. Double centralizer theorem for semisimple superalgebras
\endhead
 Now we introduce the double supercentralizer theorem for 
semisimple superalgebras. 
Although this follows easily from the discussion as in [3], 
it seems convenient to state it precisely. 
Here we partly take a finer viewpoint of strict isomorphisms. 

\proclaim{Theorem 2.1} Let $A$ be a semisimple superalgebra. 
Let $V$ be an $A$-module with the associated representation
 $\rho\:A\to\E(V)$, and 
$$
V=V_1\oplus V_2\oplus\cdots\oplus V_s\tag2.1
$$
be the decomposition of $V$ into $A$-homogeneous components. 
Put $B=\E^{\cdot}_A(V)$. We define a representation 
$\T{\rho}\:A\lot B\to\E(V)$ of $A\lot B$ by 
$\T{\rho}(a\lot b)=\rho(a)\circ b$ for all $a\in A$ and $b\in B$. 
Then each $V_i$, $1\leq i\leq s$, is a 
simple $A\lot B$-submodule of $V$ of type $M$. 

Let $U_i$ be a simple $A$-module contained in $V_i$.
 If $U_i$ is of type $M$, then there exists a simple $B$-module $W_i$ 
of type $M$ such that 
$$
U_i\lot W_i\stcong_{A\lot B}
\overline{U}_i\lot\overline{W}_i\stcong_{A\lot B}V_i.
$$

If $U_i$ is of type $Q$, then 
there exists a simple $B$-module $W_i$ of type $Q$ such that 
$$
U_i\lot W_i\stcong_{A\lot B}V_i\oplus\overline{V}_i.
$$
We have 
$$
U_i\ncong_A U_j,\quad W_i\ncong_B W_j
$$ 
for all $1\leq i\neq j\leq s$. This gives a $1-1$ correspondence between the 
isomorphism classes of simple $A$-modules and simple $B$-modules appearing 
in $V$. 
Furthermore, we have $\E^{\cdot}_B(V)=\rho(A)$. 
\endproclaim
\demo{Proof} 
Put $W_i=\H^{\cdot}_A(U_i,V_i)$ for each $i$. Since 
$b(V_i)\subset V_i$ for any $b\in B$, one can define a representation 
$\chi_i\: B\to \E^{\cdot}(W_i)$ by 
$(\chi_i(b)w)(u)=b(w(u))$ for all $b\in B$, $w\in W_i$, and $u\in U_i$. 
Then $W_i$ is a simple $B$-module. In fact, let $W_i'\subset W_i$ be 
a nonzero $B$-submodule. Since $W_i'$ is homogeneous, there is a homogeneous 
$\phi\in W_i'-\{0\}$, which can be regarded as an isomorphism onto 
$\phi(U_i)$. Let $\phi'\:\phi(U_i)\to U_i$ be its inverse. 
Also $V_i=\phi(U_i)\oplus V_i'$ for some $A$-submodule $V_i'$. 
Let $\psi$ be any homogeneous element of $W_i$. 
Define $b\in B$ by $b|_{\phi(U_i)}=\psi\circ \phi'$ and 
$b|_{V_i'}=0$, whence $\psi=b\phi\in W_i'$. Thus $W_i'$ contains all 
homogeneous elements of $W_i$, so that $W_i'=W_i$. 

If $U_i$ is of type $M$ (resp\. type $Q$) and 
$V_i\stcong_A U_i^{\oplus k_i'}\oplus \overline{U}_i^{\oplus l_i'}$ 
(resp\. $V_i\stcong_A U^{\oplus n_i'}$), then 
by Theorem 1.2 we have $\E^{\cdot}_A(V_i)\cong M(k_i',l_i')$ 
(resp\. $\E^{\cdot}_A(V_i)\cong Q(n_i')$). 
Since $\chi_i$ decomposes as 
$\E^{\cdot}_A(V)\rightarrow\E^{\cdot}_A(V_i)\rightarrow\E(W_i)$, 
the first part being a surjection and $\E^{\cdot}_A(V_i)$ being 
a simple superalgebra, we have 
$\chi_i(B)\cong \E^{\cdot}_A(V_i)$. Therefore $W_i$ is of type $M$ 
(resp\. type $Q$). 

Now $V_i$ is a quotient of $U_i\lot W_i$ by an $A\lot B$-homomorphism 
$u\lot w \to (-1)^{\alp\beta}w(u)$ for $u\in (U_i)_{\alp}$ 
and $w\in (W_i)_{\beta}$. If $U_i$ and $W_i$ are of type $M$, 
then by Theorem 1.4 $U_i\lot W_i$ is a simple $A\lot B$-module, 
so that $V_i\stcong_{A\lot B} U_i\lot W_i$. Next suppose 
$U_i$ and $W_i$ are of type $Q$. Choose 
$x_i\in \E^{1}_{A}(U_i)$ such that $x_i^2=-1$, and define 
$y_i\in \E^{1}_B(W_i)$ by
$$
(y_i(w))(u)=(-1)^{\alp}\sqrt{-1}w(x_i(u))\tag2.2
$$
for all $w\in (W_i)_{\alp}$ and $u\in U_i$. 
Then $y_i^2=-1$. By Theorem 1.4 and a comparison of 
dimensions, we have 
either $V_i\stcong_{A\lot B}(U_i\lot W_i)^{+}$ or 
$(U_i\lot W_i)^{-}$ defined with respect to $x_i\lot y_i$. 
Direct computation shows that the former is the case. 

Consequently we have 
$$
V\stcong_{A\lot B}
\bigoplus_{\t{ type } M}U_i\sot W_i
\oplus
\bigoplus_{\t{ type } Q}(U_i\sot W_i)^{+}.
\tag2.3
$$
Then, using Theorem 1.6 we have 
$$
\align
&\rho(A)\cong
\bigoplus_{\t{ type } M}M(k_i,l_i)
\oplus
\bigoplus_{\t{ type } Q}
Q(n_i)\\
\intertext{and also by Theorem 1.2} 
&B\cong 
\bigoplus_{\t{ type } M}M(k'_i,l'_i)
\oplus
\bigoplus_{\t{ type } Q}
Q(n'_i)\\
\endalign
$$
where 
$$
\alignat 3
&\vdim U_i=(k_i,l_i),\;&&\vdim W_i=(k_i',l_i')\quad &&
\t{ if $U_i$ and $W_i$ are of type $M$,}\\
&\vdim U_i=(n_i,n_i),\;&&\vdim W_i=(n_i',n_i')\quad &&
\t{ if $U_i$ and $W_i$ are of type $Q$.}\\
\endalignat
$$
In particular, the $B$-modules $W_i$ are mutually nonisomorphic. 
Therefore \thetag{2.1} can also be regarded as the decomposition of 
$V$ into $B$-homogeneous components. 
Starting with this and appropriately identifying $A\lot B$-modules 
with $B\lot A$-modules as noted in \S1, {\bf E}, we can see 
that $\E^{\cdot}_{B}(V)=\rho(A)$.
\qed\enddemo

Note that , in \thetag{2.3}, the signature in 
$(U_i\lot W_i)^{+}$ is always $+$ so long as $y_i$ is constructed 
from $x_i$ by \thetag{2.2} regardless of the choice of $x_i$. 
The above theorem gives a supercommuting action of two superalgebras 
$A$ and $B$ on $V$, which decomposes into a multiplicity-free sum of simple 
$A\lot B$-modules, where type $M$ (resp\. type $Q$) simple $A$-modules are 
paired with type $M$ (resp\. type $Q$) simple $B$-modules in a bijective 
manner. Later we also 
encounter a similar but slightly different situation, 
in which type $M$ (resp\. type 
$Q$) simple $A$-modules are paired with type $Q$ (resp\. type $M$) simple 
$B$-modules. Here we formulate this as the following corollary. 
Note that the superalgebra generated by an element $x$ of 
degree $1$ satisfying $x^2=-1$ is $2$-dimensional, and is isomorphic to 
$Q(1)$ by 
$x\mapsto 
\pmatrix
0&\sqrt{-1}\\
\sqrt{-1}&0
\endpmatrix
$. It has a unique simple module of dimension $2$, which is of type $Q$. 
(This superalgebra is also isomorphic to $\C_1$. See the 
beginning of \S3.) 
\proclaim{Corollary 2.2} Let $A$, $V$ and 
$V=V_1\oplus V_2\oplus\cdots\oplus V_s$ be as in Theorem 2.1. 
Assume that there exists $x\in\E^{1}_A(V)$ such that $x^2=-1$. 
Let $C$ denote the subsuperalgebra of $\E^{\cdot}_A(V)$ 
generated by $x$, and put 
$A'=A\lot C$ and $B=\E^{\cdot}_{A'}(V)$. 
Restricting the homomorphism 
$\rho\:A'\lot B\to \E(V)$ of Theorem 2.1, we can regard $V$ as an 
$A\lot B$-module. 
Then each $V_i$, $1\leq i\leq s$, is a 
simple $A\lot B$-submodule of $V$ of type $Q$. 

Let $T_i$ be a simple $A$-module contained in $V_i$.
 If $T_i$ is of type $M$, then there exists a simple $B$-module $W_i$ 
of type $Q$ such that 
$$
T_i\lot W_i\stcong_{A\lot B} 
\overline{T}_i\lot W_i
\stcong_{A\lot B}V_i.
$$
If $T_i$ is of type $Q$, then 
there exists a simple $B$-module $W_i$ of type $M$ such that 
$$
T_i\lot W_i\stcong_{A\lot B}
T_i\lot \overline{W}_i
\stcong_{A\lot B}V_i.
$$
We have 
$$
T_i\ncong_A T_j,\quad W_i\ncong_B W_j
$$ 
for all $1\leq i\neq j\leq s$. This also 
gives a $1-1$ correspondence between the 
isomorphism classes of simple $A$-modules and simple $B$-modules appearing 
in $V$. 

Furthermore, we have 
$\E^{\cdot}_A(V)=B\lot C$ and $\E^{\cdot}_B(V)=\rho(A)\lot C$. 
\endproclaim
\demo{Proof} From Theorem 2.1, $V$ is decomposed as a 
multiplicity-free sum of simple $A'\lot B$-modules as follows:
$$
V\stcong_{A'\lot B}
\bigoplus_{\t{ type } M}U_i\lot W_i
\oplus\bigoplus_{\t{ type } Q}
(U_i\lot W_i)^{+}.
$$
Without loss of generality, we may assume that 
the $U_i$, $1\leq i\leq m$, are of type $M$ and 
the $U_j$, $m+1\leq j\leq s$, are of type $Q$ for some $m\leq s$. 

We use Theorem 1.4 for $A$ and $C$, noting that $C$ has a unique simple 
module $X$ of dimension $2$, which is of type $Q$. 
For $1\leq i\leq m$, we are in case {\rm (c)}. 
This implies $k_i=l_i$, and that $T_i=U_i|_{A}$ 
is a simple $A$-module of type $Q$. Also by Corollary 1.5 the $T_i$, $1\leq i\leq m$, 
are all mutually nonisomorphic. For $m+1\leq i\leq s$, we are in case 
{\rm (b)}. Hence $U_i\stcong_{A'} T_i\lot X$ for some simple $A$-module 
$T_i$ of type $M$ 
with dimention $(k_i'',l_i'')$ for some $k_i''$ and $l_i''$ summing 
up to $n_i$. We can use $x_i=\id_{V_i}\lot z$, with $z\in\E^1_C(X)$ 
satisfying $z^2=-1$, for the $x_i$ required in the proof of Theorem 2.1. 
Then $(U_i\lot W_i)^+=T_i \lot 
(X \lot W_i)^+$. Since $(X \lot W_i)^+|_{B}\stcong_B W_i$, 
we have $(U_i\lot W_i)^+|_{A\lot B}\stcong_{A\lot B}T_i\lot W_i$. 
Again by Corollary 1.5, the $T_i$, $m+1\leq i\leq s$, are all mutually 
nonisomorphic. Therefore we have 
$$
V\stcong_{A\lot B}
\bigoplus_{1\leq i\leq m} 
T_i\lot W_i
\oplus
\bigoplus_{m+1\leq i\leq s} 
T_i\lot W_i. 
$$
in which all $T_i$ are distinct. Thus the $T_i\lot W_i$ are 
the $A$-homogeneous components, and this decomposition coincides with 
that into the $V_i$ in the statement. 
This establishes a $1-1$ correspondence between the simple $A$-modules 
of type $Q$ (resp\. type $M$) and the simple $B$-modules of type $M$ 
(resp\. type $Q$) appearing in $V$. 

By Theorem 2.1, we have $\E^{\cdot}_B(V)=\rho(A')$. The surjective 
homomorphism $\rho(A)\lot C\to \rho(A')$ is an isomorphism, since 
for any simple component of $S$ of $\rho(A)$, $S\lot C$ is a 
simple superalgebra by Theorem 1.3, 
and its image is present since $S\lot 1_C$ is mapped injectively. 

Moreover, clearly $\E^{\cdot}_A(V)$ contains $C$ and $B$, therefore by the 
same argument $C\lot B$. We have, by Theorem 2.1 and the above decomposition 
under $A\lot B$, 
$$
\align
B&\cong
\bigoplus_{1\leq i\leq m}M(k_i,l_i)
\oplus\bigoplus_{m+1\leq i\leq s}Q(n_i)\quad\t{ and }\\
\E^{\cdot}_A(V)&\cong
\bigoplus_{1\leq i\leq m}Q(k_i+l_i)
\oplus\bigoplus_{m+1\leq i\leq s}M(n_i,n_i).\\
\endalign
$$
Comparing the dimensions, we see that $C\lot B$ gives the whole 
$\E^{\cdot}_A(V)$. 
\qed\enddemo

\remark{Remark} The assumption in the above corollary is actually equivalent 
to the existence of an invertible element in $\E^{1}_A(V)$. In fact, 
the latter means the existence of an invertible element of degree $1$ in each 
simple component of $\E^{1}_A(V)$, which is equivalent to assuming 
$k_i'=l_i'$ for all $1\leq i\leq m$ in the decomposition of Theorem 2.1 
applied for $A$ and $\E^{\cdot}_A(V)$ (not for $A'$ and $B$). If this 
is the case, $\E^{1}_A(V)$ clearly contains an element $x'$ satisfying 
$(x')^2=1$, whence an element $x$ satisfying $x^2=-1$ also. 
The other direction is clear since such $x$ is invertible. 
\endremark 

\head
3. A simple relation between $\A_k$ and $\B_k$
\endhead
Let $\C_k$ denote the $2^k$-dimensional Clifford algebra, namely 
$\C_k$ is an associative algebra 
generated by $k$ elements $\xi_1,\dots,\xi_k$ subject to relations 
$$
\xi_i^2=1,\quad \xi_i\xi_j=-\xi_j\xi_i\quad(i\neq j)\;\;.\tag3.1
$$
We regard $\C_k$ as a superalgebra by giving degree $1$ to 
the generator $\xi_i$ for all $1\leq i\leq k$. 
For any subset 
$I=\{i_1<i_2<\cdots<i_r\}$ of $[k]=\{1,2,\dots,k\}$, we write 
$\xi_I=\xi_{i_1}\xi_{i_2}\cdots\xi_{i_r}$. Then 
$\{\xi_I\;;\;I\subset[k]\}$ is a basis of $\C_k$. 
Note that ${\Cal C}_k$ is isomorphic to 
a twisted group algebra of ${\Bbb Z}_2^k$ and 
therefore it is semisimple. Furthermore, 
it is easy to show that $\C_k$ is a simple superalgebra. 
If $k$ is even (resp\. odd), then $\C_k$ is of type $M$ (resp\. of type $Q$). 
Put $r=\lfl k/2\rfl$. 
Define a ${\Bbb Z}_2$-graded simple $\C_k$-module $X_k$ 
as a minimal left superideal (namely a minimal ${\Bbb Z}_2$-graded left ideal) 
of $\C_k$ as follows:
$$
\align
&X_k=\C_ke,\tag3.2\\ 
&e=e_1e_2\cdots e_r,\quad 
e_i=\frac{1}{\sqrt{2}}(1+\sqrt{-1}\xi_{2i-1}\xi_{2i}).\\
\endalign
$$
For each $\vep=(\vep_1,\dots,\vep_r)\in{\Bbb Z}_2^r$, put 
$\xi^{\vep}=(\xi_1^{\vep_1}e_1)(\xi_3^{\vep_2}e_2)\cdots
(\xi_{2r-1}^{\vep_r}e_r)$. If $k$ is even (resp\. odd), 
then $\{\xi^{\vep}\,;\,\vep\in{\Bbb Z}_2^r\}$ 
(resp\. $\{\xi^{\vep}\,;\,\vep\in{\Bbb Z}_2^r\}
\cup\{\xi^{\vep}\xi_k\,;\,\vep\in{\Bbb Z}_2^r\}$) 
is a basis of $X_k$. 
If $k$ is odd, define $z_k\in\E^{1}_{\C_k}(X_k)$ by 
$$
z_k(\xi^{\vep}\xi_k^{\alp})
=(-1)^{\sum_i\vep_i+\alp}(\xi^{\vep}\xi_k^{\alp+1})
$$
for all $\vep=(\vep_1,\dots,\vep_r)\in{\Bbb Z}_2^r$ and 
$\alp\in{\Bbb Z}_2$. Note that $z_k^2=-1$.

\proclaim{Proposition 3.1} \t{(cf\. [12, Prop\. 3.1])} The character value of 
the simple $\C_k$-module $X_k$ is given by 
$$
\Ch[X_k]\left(\sum_Ic_I\xi_I\right)=2^{\lfloor (k+1)/2\rfloor}c_{\phi}
$$
for any element $\sum_Ic_I\xi_I\in\C_k$, with $c_I\in{\Bbb C}$ for all 
$I\subset[k]$. 
\endproclaim

The following isomorphism is the key to our clarification 
of the relation between the simple modules over $\A_k$ and $\B_k$. 
\proclaim{Theorem 3.2} There exists a superalgebra homomorphism 
$\vth\:\C_k\lot\A_k\to\B_k$ satisfying
$$
\align
&\vth(\xi_1\lot1)\mapsto \tau_i\;\;\;(1\leq i\leq k),\tag3.3\\
&\vth(1\lot\gam_j)\mapsto
\frac{1}{\sqrt{2}}(\tau_j-\tau_{j+1})\sig_j
\;\;\;(1\leq j\leq k-1)
\endalign
$$
where $\tau_i=\sig_{i-1}\cdots\sig_1\tau\sig_1\dots\sig_{i-1}$ 
for all $1\leq i\leq k$. Then it follows that $\vth$ is an isomorphism.
\endproclaim
\demo{Proof} Put 
$$
\T{\gam}_j=
\dfrac{1}{\sqrt{2}}(\tau_j-\tau_{j+1})\sig_j
$$ 
for all $1\leq j\leq k-1$. 

It is easy to check that the $\tau_i$, $1\leq i\leq k$, and 
the $\T{\gam}_j$, $1\leq j\leq k-1$, satisfy 
the relations \thetag{3.1} and \thetag{1.2} in which $\xi_i$ and $\gam_j$ 
are replaced by $\tau_i$ and $\T{\gam}_j$ respectively. 
Furthermore, we have $\tau_i\T{\gam}_j=-\T{\gam}_j\tau_i$ 
for all $1\leq i\leq k$, $1\leq j\leq k-1$, since 
$$
\align
&\tau_i
(\tau_j-\tau_{j+1})=\left\{
\matrix
\format\l&\quad\l\\
-(\tau_j-\tau_{j+1})\tau_{i+1}&\t{if $i=j$,}\\
-(\tau_j-\tau_{j+1})\tau_{i-1}&\t{if $i=j+1$,}\\
-(\tau_j-\tau_{j+1})\tau_i&\t{if $i\neq j,j+1$,}
\endmatrix\right.\\
&\tau_i\sig_j=\left\{
\matrix
\format\l&\quad\l\\
\sig_j\tau_{i+1}&\t{if $i=j$,}\\
\sig_j\tau_{i-1}&\t{if $i=j+1$,}\\
\sig_j\tau_i&\t{if $i\neq j,j+1$.}
\endmatrix\right.\\
\endalign
$$ 
Therefore $\vth$ is a homomorphism. 

Since 
$$
\vth(\xi_1\lot1)=\tau\quad\t{and}\quad
\vth(\frac{1}{\sqrt{2}}(\xi_j-\xi_{j+1})\lot\gam_j)=\sig_j
$$
for all $1\leq j\leq k-1$, it follows that $\vth$ 
is surjective. Furthermore, it follows that $\vth$ is bijective, 
since $\dim{(\C_k\lot\A_k)}=\dim{\B_k}=2^kk!$. 
Therefore $\vth$ is an isomorphism. \qed
\enddemo

Using the above isomorphism, 
we identify $\C_k\lot\A_k$ with $\B_k$ and 
regard $\C_k\lot\A_k$-modules $X_k\lot V_{\nu}$, 
$\nu\in DP_k$, as $\B_k$-modules. 
By Theorem 1.4, we can construct all simple $\B_k$-modules as follows. 
For each $\nu\in DP_k^{-}$, fix a nonzero element $x_{\nu}$ of 
$\E^{1}_{\A_k}(V_{\nu})$. 

\proclaim{Proposition 3.3} {\rm(1)} If $k$ is even, then 
$\{X_k\lot V_{\nu},$ $X_k\lot\overline{V}_{\nu}\,|\,\nu\in DP_k^{+}\}$ 
(resp\. 
$\{X_k\lot V_{\nu}\,|\,\nu\in DP_k^{-}\}$) is a complete set 
of strict isomorphism classes of simple $\B_k$-modules of type $M$ 
(resp\. type $Q$). 

{\rm(2)} If $k$ is odd, then 
$\{(X_k\lot V_{\nu})^{+},$ $(X_k\lot V_{\nu})^{-}\,|\,
\nu\in DP_k^{-}\}$ 
(resp\. $\{X_k\lot V_{\nu}\,|\,\nu\in DP_k^{+}\}$) is a complete 
set of strict isomorphism clases of simple $\B_k$-modules of type $M$ 
(resp\. type $Q$), where, if $\nu\in DP_k^{-}$, 
$(X_k\lot V_{\nu})^{\pm}$ denotes the $\pm\sqrt{-1}$-eigenspace of 
$z_k\lot x_{\nu}\in\E^{0}_{\B_k}(X_k\lot V_{\nu})$ respectively. 
\endproclaim

\comment
It is easy to check that 
$$
\overline{(X_k\lot V_{\nu})^{+}}\stcong_{\B_k}
(X_k\lot V_{\nu})^{-}.
$$
\endcomment
 
The following proposition shows that 
our parametrization of the simple $\B_k$-modules 
coincides with Sergeev's parametrization in Theorem 1.9. 

\proclaim{Proposition 3.4} If $k$ is even or $\nu\in DP_k^{+}$, then 
$$
\Ch[X_k\lot V_{\nu}]=\psi_{\nu}.
$$
If $k$ is odd and $\nu\in DP_k^{-}$, then 
$$
\Ch[(X_k\lot V_{\nu})^{+}]
=\Ch[(X_k\lot V_{\nu})^{-}]=\psi_{\nu}.
$$ 
Here $\psi_{\nu}$ is the irreducible character defined in Theorem 1.9.  
\endproclaim
\demo{Proof} 
Let $\nu\in DP_k$, $\mu\in OP_k$ and put $l=l(\mu)$. Since 
$\sig^{(\mu,\phi)}$ is a product of $k-l$ generators $\sig_j$, its image 
under $\vth^{-1}$ is a product of $k-l$ elements of $\C_k\lot\A_k$ 
of the form $\frac{1}{\sqrt{2}}(\xi_j-\xi_{j+1})\lot \gam_j$. Rearrange 
this product into the form 
$$
\t{(constant)}\times
\t{(product of the $\xi_j-\xi_{j+1}$)}\lot
\t{(product of the $\gam_j$)}.
$$
Let $\vX$ denote the product of the $\xi_j-\xi_{j+1}$ in this expression. 
The product of the $\gam_j$ equals $\gam^{\mu}$. 
Therefore we have 
$$
\vth^{-1}(\sig^{(\mu,\phi)})
=(-1)^{\binom{k-l}{2}}\left(\frac{1}{\sqrt{2}}\right)^{k-l}\vX\lot\gam^{\mu}, 
$$
where the signature comes from interchanging the elements of degree $1$. 
This gives 
$$
\Ch[X_k\lot V_{\nu}](\vth^{-1}
(\sig^{(\mu,\phi)})=
(-1)^{\binom{k-l}{2}}
\left(\frac{1}{\sqrt{2}}\right)^{k-l}
\Ch[X_k](\vX)\Ch[V_{\nu}](\gam^{\mu}). 
$$
By Proposition 3.1, $\Ch[X_k](\vX)$ equals 
$2^{\lfl(k+1)/2\rfl}$ 
times the coefficient of $1$ in the expansion of this product into the 
$\xi_I$, $I\subset[k]$. We have 
$$
\vX=\alp_1(\vX_1)\alp_2(\vX_2)\cdots\alp_l(\vX_l)
$$
where 
$$
\vX_j=(\xi_1-\xi_2)(\xi_2-\xi_3)\cdots(\xi_{\mu_j-1}-\xi_{\mu_j})\in\C_{\mu_j}
$$
and $\alp_j\:\C_{\mu_j}\to\C_k$ is the embedding 
defined by $\xi_i\mapsto\xi_{\mu_1+\mu_2+\cdots+\mu_{j-1}+i}$ 
$(1\leq i\leq \mu_j, 1\leq j\leq l)$. 
If we let $c_j$ denote the coefficient of $1$ in the expansion 
of $\vX_j$ into the $\xi_{I_j}$, $I_j\subset[\mu_j]$, then 
we have $\Ch[X_k](\vX)=2^{\lfl(k+1)/2\rfl}c_1c_2\cdots c_l$. Clearly 
$$
c_j=\left\{
\matrix
\format\l&\quad\l\\
0
&\t{if $\mu_j$ is even,}\\
(-1)^{\frac{\mu_j-1}{2}}
&\t{if $\mu_j$ is odd,}
\endmatrix\right.
$$
so that we have 
$\Ch[X_k](\vX)=2^{\lfl(k+1)/2\rfl}(-1)^{(k-l)/2}$. 
This signature cancels with 
$(-1)^{\binom{k-l}{2}}$, since 
$\binom{k-l}{2}+\frac{k-l}{2}=\frac{(k-l)^2}{2}$ and we have 
$k-l\equiv0$ (mod $2$) (so that $(k-l)^2\equiv0$ (mod $4$)) because 
$\mu\in OP_k$. Therefore we have 
$$
\align
&\Ch[X_k\lot V_{\nu}]
(\vth^{-1}(\sig^{(\mu,\phi)}))\\
&\quad=
2^{\lfl\frac{k+1}{2}\rfl}(\sqrt{2})^{l-k}
\Ch[V_{\nu}](\gam^{\mu})
=\left\{
\matrix
\format\l&\quad\l\\
(\sqrt{2})^{l}\Ch[V_{\nu}](\gam^{\mu})&\t{if $k$ is even,}\\
&\\
(\sqrt{2})^{l+1}\Ch[V_{\nu}](\gam^{\mu})&
\t{if $k$ is odd.}\\
\endmatrix\right.
\endalign
$$
In particular, if $k$ is odd and $\nu\in DP_k^{-}$, we have 
$$
\align
&\Ch[(X_k\lot V_{\nu})^{+}]
(\vth^{-1}(\sig^{(\mu,\phi)}))
=\Ch[(X_k\lot V_{\nu})^{-}]
(\vth^{-1}(\sig^{(\mu,\phi)}))\\
&=(\sqrt{2})^{l-1}\Ch[V_{\nu}](\gam^{\mu}).
\endalign
$$

If $k$ is even, then $d(\nu)=\vep(\nu)$ for any $\nu\in DP_k$. 
By Theorem 1.8, Theorem 1.9 and Theorem 1.10, we have 
$$
\align
&\sum_{\nu\in DP_k}
\Ch[X_k\lot V_{\nu}]
(\vth^{-1}(\sig^{(\mu,\phi)}))
(\sqrt{2})^{-l(\nu)-d(\nu)}Q_{\nu}\\
&\quad=(\sqrt{2})^{l(\mu)}
\sum_{\nu\in DP_k}\Ch[V_{\nu}](\gamma^\mu)
(\sqrt{2})^{-l(\nu)-\vep(\nu)}Q_{\nu}\\
&\quad=2^{l(\mu)}p_{\mu}\\
&\quad=\sum_{\nu\in DP_k}\psi_{\nu}(\sig^{(\mu,\phi)})
(\sqrt{2})^{-l(\nu)-d(\nu)}Q_{\nu}.
\endalign
$$
Since the $Q_{\nu}$, $\nu\in DP_k$, are linearly 
independent in $\Omega^k$ (cf\. \S1, {\bf C}), we have  
$\Ch[X_k\lot V_{\nu}](\vth^{-1}
(\sig^{(\mu,\phi)}))=\psi_{\nu}(\sig^{(\mu,\phi)})$.
 Therefore we have $\Ch[X_k\lot V_{\nu}]=\psi_{\nu}$. 

If $k$ is odd, then $d(\nu)+\vep(\nu)=1$ for any $\nu\in DP_k$. 
We have 
$$
\align
&\sum_{\nu\in DP_k}
2^{-\vep(\nu)}\Ch[X_k\lot V_{\nu}]
(\vth^{-1}(\sig^{(\mu,\phi)}))
(\sqrt{2})^{-l(\nu)-d(\nu)}Q_{\nu}\\
&=\sum_{\nu\in DP_k}
(\sqrt{2})^{l+1-2\vep(\nu)}
\Ch[V_{\nu}](\gam^{\mu})
(\sqrt{2})^{-l(\nu)-d(\nu)}Q_{\nu}\\
&\quad=2^{l(\mu)}p_{\mu}\\
&\quad=\sum_{\nu\in DP_k}\psi_{\nu}(\sig^{(\mu,\phi)})
(\sqrt{2})^{-l(\nu)-d(\nu)}Q_{\nu}.
\endalign
$$
This gives $2^{-\vep(\nu)}\Ch[X_k\lot V_{\nu}]=\psi_{\nu}$. 
Namely, $\Ch[X_k\lot V_{\nu}]=\psi_{\nu}$ if $\nu\in DP_k^{+}$ and 
$\Ch[(X_k\lot V_{\nu})^{+}]=
\Ch[(X_k\lot V_{\nu})^{-}]
=\psi_{\nu}$ if $\nu\in DP_k^{-}$. 
\qed\enddemo

By Proposition 3.4 , we can fix the choice of simple 
$\B_k$-modules $W_{\nu}$, 
which afford the characters $\psi_{\nu}$, as follows: 
$$
W_{\nu}=\left\{
\matrix
\format
\l&\quad\l\\
X_k\lot V_{\nu}&\t{if $k$ is even or $\nu\in DP_k^{+}$,}\\
(X_k\lot V_{\nu})^{+}&\t{if $k$ is odd and $\nu\in DP_k^{-}$,}
\endmatrix\right.\tag3.4
$$
where $(X_k\sot V_{\nu})^{\pm}$ are as in Proposition 3.3. 
The $W_{\nu}$, $\nu\in DP_k$, form a complete set of (not strict) 
isomorphism classes of simple $\B_k$-module. 

The following proposition gives 
the restriction rule of the simple $\B_k$-modules to 
$\A_k$ by $\A_k\overset{\cong}\to{\rightarrow}
1\lot\A_k\subset\C_k\lot\A_k\overset{\vth}\to{\rightarrow}\B_k$. 

\proclaim{Proposition 3.5} For any $\nu\in DP_k$, 
the simple $\B_k$-module $W_{\nu}$ restricts to an $\A_k$-module 
as follows:
$$
W_{\nu}|_{\A_k}
\stcong
\overline{W_{\nu}|_{\A_k}}
\stcong
\left\{
\matrix
\format\l&\l\\
(V_{\nu}\oplus\overline{V}_{\nu})^{\oplus 2^{\lfl(k-1)/2\rfl}}\;&
\;\t{\rm if $\nu\in DP^{+}_k$,}\\ 
V_{\nu}^{\oplus 2^{\lfl (k+1)/2\rfl}}\;&
\;\t{\rm if $\nu\in DP^{-}_k$, $k$ is even,}\\ 
V_{\nu}^{\oplus 2^{\lfl (k-1)/2\rfl}}\;&
\;\t{\rm if $\nu\in DP^{-}_k$, $k$ is odd,}\\ 
\endmatrix
\right.
$$
where 
$\overline{W_{\nu}|_{\A_k}}$ and 
$\overline{V}_{\nu}$ are shifts of $W_{\nu}|_{\A_k}$ and 
$V_{\nu}$ respectively. 
\endproclaim
\demo{Proof} First, in view of Theorem 1.8, 
the $\A_k$-modules on the right-hand sides are all strictly isomorphic to 
their shifts. 
Assume $\nu\in DP^{+}_k$. We have 
$\xi\lot V_{\nu}\stcong V_{\nu}$ (resp\. $\overline{V}_{\nu}$) 
if $\xi\in(X_k)_{0}$ (resp\. $(X_k)_{1}$). 
From \thetag{3.4} and the fact that 
$\dim{(X_k)_0}=\dim{(X_k)_1}=2^{\lfl (k-1)/2\rfl}$, we have 
$$
W_{\nu}|_{\A_k}\stcong V_{\nu}^{\oplus 2^{\lfl(k-1)/2\rfl}} 
\oplus \overline{V}_{\nu}^{\oplus 2^{\lfl(k-1)/2\rfl}}.
$$

Assume $\nu\in DP_k^{-}$. We have 
$\xi\lot V_{\nu}\stcong V_{\nu}\stcong\overline{V}_{\nu}$ 
for any homogeneous element $\xi\in X_k$. 
Therefore we have 
$W_{\nu}|_{\A_k}=V_{\nu}^{\oplus 2^{\lfl (k+1)/2\rfl}}$ if $k$ is even, and 
we have $W_{\nu}|_{\A_k}\oplus\overline{W}_{\nu}|_{\A_k}
=X_k\lot V_{\nu}|_{\A_k}
=V_{\nu}^{\oplus 2^{\lfl (k+1)/2\rfl}}$ if $k$ is odd. 
Note that $\overline{W}_{\nu}=(X_k\lot V_{\nu})^{-}$. 
The result also follows in this case. 
\qed\enddemo

\head 
4. A duality relation of $\A_k$ and $\fq(n)$
\endhead
In this section, we establish 
a duality relation between $\A_k$ and $\U_n$. 

We can rewrite \thetag{1.11} using the description of 
$W_{\nu}$ after the Proof of Proposition 3.4. If 
$k$ is even, then we have 
$$
W\stcong
\bigoplus_{\nu\in DP_k^{+}}
(X_k\lot V_{\nu})\lot U_{\nu}\oplus
\bigoplus_{\nu\in DP_k^{-}}
\left((X_k\lot V_{\nu})\lot U_{\nu}\right)^{+}
$$
and if $k$ is odd, then we have 
$$
W\stcong
\bigoplus_{\nu\in DP_k^{-}}
(X_k\lot V_{\nu})^{+}\lot U_{\nu}
\oplus
\bigoplus_{\nu\in DP_k^{+}}
\left((X_k\lot V_{\nu})\lot U_{\nu}\right)^{+}.
$$
Note that the symbol $^{+}$ is used in three ways. 
\roster
\item If $k$ is odd and $\nu\in DP_k^{-}$, then both 
$X_k$ and $V_{\nu}$ are of type $Q$. 
$(X_k\lot V_{\nu})^{+}$ is the $(+\sqrt{-1})$-eigenspace of 
$z_k\lot x_{\nu}$. 
\item If $k$ is even and $\nu\in DP_k^{-}$, then 
$X_k$, $V_{\nu}$, $U_{\nu}$ are of type $M$, $Q$, $Q$ respectively. 
$((X_k\lot V_{\nu})\lot U_{\nu})^{+}$ 
is the $(+\sqrt{-1})$-eigenspace of 
$(1\lot x_{\nu})\lot u_{\nu}$, where 
$u_{\nu}\in\E^1_{\U_n}(U_{\nu})$ is defined from $1\lot x_{\nu}$ as in 
\thetag{2.2}.
\item If $k$ is odd and $\nu\in DP_k^{+}$, then 
$X_k$, $V_{\nu}$, $U_{\nu}$ are of type $Q$, $M$, $Q$ respectively. 
$((X_k\lot V_{\nu})\lot U_{\nu})^{+}$ 
is the $(+\sqrt{-1})$-eigenspace of 
$(z_k\lot 1)\lot u_{\nu}$, where 
$u_{\nu}\in\E^1_{\U_n}(U_{\nu})$ is defined from $z_k\lot1$ as in 
\thetag{2.2}.
\endroster

Put $r=\lfl k/2\rfl$ and 
$\zeta_i=\sqrt{-1}\xi_{2i-1}\xi_{2i}\in\C_k$ for 
$1\leq i\leq r$. The $\Psi(\zeta_i)$, $1\leq i\leq r$, are 
commuting involutions of 
$\Psi((\C_k)_0)\subset\Psi((\B_k)_0)=\E^{0}_{\Theta(\U_n)}(W)$. 

Then $W$ is a direct sum of the simultaneous eigenspaces 
$W^{\vep}$, $\vep=(\vep_1,\dots,\vep_r)
\in{\Bbb Z}_2^r$, namely we have 
$$
\align
&W=\bigoplus_{\vep=(\vep_1,\dots,\vep_r)\in{\Bbb Z}_2^r}W^{\vep},\\
&W^{\vep}=\{w\in W\,;\,\Psi(\zeta_i)(w)=(-1)^{\vep_i}w
\quad(1\leq i\leq r)\}.
\endalign
$$
It is clear that $W^{\vep}$ is an $\A_k\lot\U_n$-module for each 
$\vep\in{\Bbb Z}_2^r$. 

\proclaim{Theorem 4.1} 
For each $\vep\in{\Bbb Z}_2^r$, the submodule $W^{\vep}$ 
is decomposed as a multiplicity-free sum of simple 
$\A_k\lot\U_n$-modules, in which the simple $\A_k$-modules are 
paired with the simple $\U_n$-modules in a bijective manner 
(\thetag{4.1} and \thetag{4.3}). 
More precisely, we have the following decomposition. 

For each $\vep=(\vep_1,\dots,\vep_r)\in{\Bbb Z}_2^r$, 
put $\alp(\vep)=\vep_1+\cdots+\vep_r\in{\Bbb Z}_2$. 

\t{\rm(1)} Assume that $k$ is even. 
If $\alp(\vep)=0$, then $W^{\vep}$ is a direct sum of simple 
$\A_k\lot\U_n$-modules of type $M$ as follows: 
$$
W^{\vep}\stcong_{\A_k\lot\U_n}
\bigoplus_{\nu\in DP_k^{+}}
V_{\nu}\sot U_{\nu}\oplus
\bigoplus_{\nu\in DP_k^{-}}
(V_{\nu}\sot U_{\nu})^{+}
\tag4.1
$$
where $(V_{\nu}\sot U_{\nu})^{\pm}$ 
denotes the $\pm\sqrt{-1}$-eigenspace of 
$x_{\nu}\lot u_{\nu}\in\E^{0}_{\A_k\lot\U_n}(V_{\nu}\lot U_{\nu})$ 
, where $u_{\nu}$ is defined as in Theorem 1.10 with 
$y_{\nu}=1\lot x_{\nu}$ (if $k$ is even and $\nu\in DP_k^{-}$) or 
$y_{\nu}=z_k\lot1$ (if $k$ is odd and $\nu\in DP_k^{+}$). 

If $\alp(\vep)=1$, then $W^{\vep}$ is also a direct sum of simple 
$\A_k\lot\U_n$-modules of type $M$ and 
is strictly isomorphic to the shift 
of the module on the right-hand side in \thetag{4.1}. 

Furthermore we have 
$$
\E^{\cdot}_{\Theta(U_n)}(W^{\vep})=\Psi(\A_k),\quad\;
\E^{\cdot}_{\Psi(\A_k)}(W^{\vep})=\Theta(\U_n)\tag4.2
$$
for all $\vep\in{\Bbb Z}_2^r$.

\t{\rm(2)} Assume that $k$ is odd. 
Then $W^{\vep}$ is an $\A_k\lot\U_n$-module of type $Q$ and 
$$
W^{\vep}\cong_{\A_k\lot\U_n}
\bigoplus_{\nu\in DP_k}V_{\nu}\sot U_{\nu}. 
\tag4.3
$$

Furthermore we have 
$$
\E^{\cdot}_{\Theta(U_n)}(W^{\vep})\cong\C_1\otimes\Psi(\A_k),\quad
\E^{\cdot}_{\Psi(\A_k)}(W^{\vep})\cong\C_1\otimes\Theta(\U_n).\tag4.4 
$$
\endproclaim
\demo{Proof} 
For each $\vep\in{\Bbb Z}_2^r$, let $X_k^{\vep}$ denote 
the simultaneous eigenspace of $X_k$ of the 
$\zeta_i=\sqrt{-1}\xi_{2i-1}\xi_{2i}$, $1\leq i\leq r$, namely 
$$
X_k^{\vep}=
\{\xi\in X_k\,;\,\zeta_i\xi=(-1)^{\vep_i}\xi
\quad(1\leq i\leq r)\}.
$$
(1) Assume that $k$ is even. 
Then we have $X_k^{\vep}={\Bbb C}\xi^{\vep}$ 
since $\zeta_i\xi^{\vep}=(-1)^{\vep_i}\xi^{\vep}$ 
for each $\vep\in{\Bbb Z}_2^r$. 
\comment
$$
\align
\zeta_i\xi^{\vep}
&=(\xi_1^{\vep_1}e_1)\cdots(\xi_{2i-3}^{\vep_{i-1}}e_{i-1})
(\zeta_i\xi_{2i-1}^{\vep_i}e_i)
(\xi_{2i+1}^{\vep_{i+1}}e_{i+1})\cdots
(\xi_{2r-1}^{\vep_r}e_r)\\
&=(-1)^{\vep_i}(\xi_1^{\vep_1}e_1)
\cdots(\xi_{2i-3}^{\vep_{i-1}}e_{i-1})
(\xi_{2i-1}^{\vep_i}e_i)
(\xi_{2i+1}^{\vep_{i+1}}e_{i+1})\cdots
(\xi_{2r-1}^{\vep_r}e_r)
=(-1)^{\vep_i}\xi^{\vep}
\endalign
$$
\endcomment

If $\nu\in DP_k^{-}$, then we have 
$$
(X_k\lot V_{\nu}\lot U_{\nu})^{+}
=\bigoplus_{\vep\in{\Bbb Z}_2^r}
X_k^{\vep}\lot(V_{\nu}\lot U_{\nu})^{+}
$$
since the $\zeta_i$, $1\leq i\leq r$, 
and $1\lot x_{\nu}\lot u_{\nu}$ commute. 
If $\nu\in DP_k^{+}$, then we have 
$$
X_k\lot V_{\nu}\lot U_{\nu}
=\bigoplus_{\vep\in{\Bbb Z}_2^r}
X_k^{\vep}\lot(V_{\nu}\lot U_{\nu}).
$$
Consequently, 
$$
\align
W^{\vep}&
\stcong_{\A_k\lot\U_n}
\xi^{\vep}\lot
\left(\bigoplus_{\nu\in DP_k^{+}}
V_{\nu}\lot U_{\nu}
\oplus\bigoplus_{\nu\in DP_k^{-}}
(V_{\nu}\lot U_{\nu})^{+}\right)\\
&\stcong_{\A_k\lot\U_n}
\left\{
\matrix
\format\l&\quad\l\\
\bigoplus_{\nu\in DP_k^{+}}
V_{\nu}\lot U_{\nu}
\oplus\bigoplus_{\nu\in DP_k^{-}}
(V_{\nu}\lot U_{\nu})^{+}&\t{if $\alp(\vep)=0$,}\\
&\\
\overline{\bigoplus_{\nu\in DP_k^{+}}
V_{\nu}\lot U_{\nu}
\oplus\bigoplus_{\nu\in DP_k^{-}}
(V_{\nu}\lot U_{\nu})^{+}}&\t{if $\alp(\vep)=1$.}\\
\endmatrix\right.
\endalign
$$
Therefore \thetag{4.1} follows. 

Since $W^{\vep}$ is an $\A_k\lot\U_n$-module, we have 
$$
\Theta(\U_n)|_{W^{\vep}}\subset\E^{\cdot}_{\Psi(\A_k)}(W^{\vep}).
$$ 
Define a linear map 
$\pr_{\vep}\:\Theta(\U_n)\to\Theta(\U_n)|_{W^{\vep}}$ by 
$\pr_{\vep}(f)=f|_{W^{\vep}}$ for all $f\in\Theta(\U_n)$. 
We claim that $\pr_{\vep}$ is injective. 
Assume that $f\in\ker\pr_{\vep}$, namely $\pr_{\vep}(f)=0\in\E(W^{\vep})$. 
Since $f$ and $\xi_{2j-1}$'s commute, and 
a subgroup of $(\C_k)^{\times}$ generated by 
the $\xi_{2j-1}$, $1\leq j\leq r$, transitively act 
on $\{W^{\vep'}\,;\,\vep'\in{\Bbb Z}_2^r\}$ as follows:
$$
\xi_{2j-1}W^{(\vep_1,\dots,\vep_r)}
=W^{(\vep_1,\dots,\vep_j+1,\dots,\vep_r)}
$$
for all $1\leq j\leq r$, 
it follows that $f|_{W^{\vep'}}=0$ for all $\vep'\in{\Bbb Z}_2^r$. 
Therefore $f=0$ in $\E(W)$, as required. 
Hence $\pr_{\vep}$ is injective. 

From Theorem 1.10, Theorem 2.1, and \thetag{4.1}, we have 
$$
\dim\E^{\cdot}_{\Psi(\A_k)}(W^{\vep})=
\dim\E^{\cdot}_{\Psi(\B_k)}(W)
=\dim\Theta(\U_n)=\dim\Theta(\U_n)|_{W^{\vep}}.
$$
Consequently we have 
$$
\E^{\cdot}_{\Psi(\A_k)}(W^{\vep})=\Theta(\U_n)|_{W^{\vep}}
$$
and, from Theorem 2.1, we have
$$
\E^{\cdot}_{\Theta(\U_n)}(W^{\vep})=\Psi(\A_k).
$$

(2) Assume that $k$ is odd. 
If $\nu\in DP_k^{+}$, then we 
regard the $\A_k\lot\C_k$-module $V_{\nu}\lot X_k$ 
as a $\C_k\lot\A_k$-module via $\omega_{\C_k,\A_k}$ (cf. \S1, {\bf E}). 
Then $X_k\lot V_{\nu}\stcong V_{\nu}\lot X_k$ 
by the map 
$\theta\:X_k\lot V_{\nu}\to V_{\nu}\lot X_k$ defined by 
$\theta(\xi\lot v)=(-1)^{\alp\cdot\beta}v\lot \xi$ for 
all homogeneous $\xi\in (X_k)_{\alp}$ 
and $v\in (V_{\nu})_{\beta}$ ($\alp$, $\beta\in{\Bbb Z}_2$). 
Since $\theta\circ(z_k\lot1)=(1\lot z_k)\circ\theta$, we have 
$$
(X_k\lot V_{\nu}\lot U_{\nu})^{+}
\stcong_{\B_k\lot\U_n}
V_{\nu}\lot(X_k\lot U_{\nu})^{+}
$$
where $(X_k\lot U_{\nu})^{\pm}$ 
denotes the $\pm\sqrt{-1}$-eigenspace of 
$z_k\lot u_{\nu}\in\E^{0}_{\C_k\lot\U_n}(X_k\lot U_{\nu})$ 
respectively. 
Since the $\zeta_i$, $1\leq i\leq r$, 
and $z_k\lot u_{\nu}$ commute, we have 
$$
(X_k\lot U_{\nu})^{+}=
\bigoplus_{\vep\in{\Bbb Z}_2^r}
(X_k^{\vep}\lot U_{\nu})^{+}
$$
where $(X_k^{\vep}\lot U_{\nu})^{\pm}$ 
denotes the $\pm\sqrt{-1}$-eigenspace of 
$z_k|_{X_k^{\vep}}\lot u_{\nu}\in
\E^{0}_{\U_n}(X_k^{\vep}\lot U_{\nu})$ respectively. 
Note that we have 
$X_k^{\vep}={\Bbb C}\xi^{\vep}\oplus{\Bbb C}\xi^{\vep}\xi_k$, 
since $\zeta_i(\xi^{\vep}\xi_k^{\alp})
=(-1)^{\vep_i}\xi^{\vep}\xi_k^{\alp}$ 
for each $\vep\in{\Bbb Z}_2^r$ and $\alp\in{\Bbb Z}_2$. 
Then it is clear that 
$(X_k^{\vep}\lot U_{\nu})^{+}$ is a 
$\U_n$-module for each $\vep\in{\Bbb Z}_2^r$. 
Moreover, by theorem 1.4 (c) we have 
$$
(X_k\lot U_{\nu})^{+}\stcong_{\U_n}U_{\nu}^{\oplus 2^r}. 
$$
Therefore we have  
$$
(X_k^{\vep}\lot U_{\nu})^{+}\stcong_{\U_n} U_{\nu}
$$
for all $\vep\in{\Bbb Z}_2^r$. 
If $\nu\in DP_k^{-}$, then we have 
$$
(X_k\lot V_{\nu})^{+}\lot U_{\nu}
=\bigoplus_{\vep\in{\Bbb Z}_2^r}
(X_k^{\vep}\lot V_{\nu})^{+}\lot U_{\nu}
$$
since the $\zeta_i$, $1\leq i\leq r$, and $(z_k\lot x_{\nu})\lot1$ commute, 
where $(X_k^{\vep}\lot V_{\nu})^{\pm}$ 
denotes the $\pm\sqrt{-1}$-eigenspace of 
$z_k|_{X_k^{\vep}}\lot x_{\nu}
\in\E^{0}_{\A_k}(X_k^{\vep}\lot V_{\nu})$ respectively. 
By Theorem 1.4 (c) we have 
$$
(X_k^{\vep}\lot V_{\nu})^{+}\stcong_{\A_k} V_{\nu}.
$$
Consequently we have 
$$
W^{\vep}
\stcong
\overline{W^{\vep}}
\stcong
\bigoplus_{\nu\in DP_k}
V_{\nu}\lot U_{\nu}.
$$
Therefore \thetag{4.3} 
follows. From Corollary 2.2, \thetag{4,4} follows. \qed
\enddemo
\proclaim{Corollary 4.2} The duality relation 
of $\A_k$ and $\U_n$ in Theorem 4.1 gives Schur's formula \thetag{1.6}. 
Namely we get Schur's formula \thetag{1.6} 
by the calculation of the character 
of the $\A_k\lot\U_n$-module $W^{\vep}$ in Theorem 4.1. 
\endproclaim
\demo{Proof} 
By what we noted before Theorem 1.10, any $\A_k\lot\U_n$-submodule $W'$ 
of $W$ can be regarded as an $\A_k$-module with a commuting polynomial 
representation $\theta_{W'}$ of $GL(n,{\Bbb C})$. Here we extend our 
notation in Theorem 1.10 to let 
$\Ch[W']\in Z_0((\A_k\otimes{\Bbb C}[GL(n,{\Bbb C})]^{*})$ be determined by 
$x\otimes g\mapsto \tr(x_{W'}\circ\theta_{W'}(g))$ for $x\in\A_k$ and 
$g\in GL(n,{\Bbb C})$, where $x_{W'}$ denotes the action of $x\in\A_k$ 
on $W'$. 

For any $\vep$, $\vep'\in{\Bbb Z}_2^r$, 
we have $\Ch[W^{\vep}]=\Ch[W^{\vep'}]$, 
since $W^{\vep}\cong W^{\vep'}$. 
Then, we have $\Ch[W]=2^r\Ch[W^{\vep}]$ 
for any $\vep\in{\Bbb Z}_2^r$. 
Therefore, for each $\mu\in OP_k$ and each diagonal element 
$E=\t{\rm diag\,}(x_1,x_2,\dots,x_n)\in GL(n,{\Bbb C})$, 
we have 
$$
\Ch[W^{\vep}]\left(\gam^{\mu}\otimes E\right)
=2^{-r}\Ch[W]\left(\vth(1\lot\gam^{\mu})\otimes E\right).
$$
Put $l=l(\mu)$. Since $1\lot\gam^{\mu}$ is a product of $k-l$ elements 
$1\lot\gam_j$, its image under $\vth$ is a product of $k-l$ 
elements of $\B_k$ of the form $\frac{1}{\sqrt{2}}(\tau_j-\tau_{j+1})\sig_j$. 
Rearrange this product into the form 
$$
(\t{constant})
\times(\t{product of the $\tau_{r}-\tau_{s}$})\times 
(\t{product of the $\sig_j$}).
$$
The product of the $\sig_j$ equals $\sig^{(\mu,\phi)}$. 
Expanding the product of $\tau_r-\tau_s$ 
into a sum of $k-l$ elements, we have 
$$
\vth(1\lot\gam^{\mu})
=\left(\frac{1}{\sqrt{2}}\right)^{k-l}\times \sum 
(\t{product of the $\tau_{r}$})\times \sig^{(\mu,\phi)}.
$$
Then all terms in the summation are conjugate to 
$\sig^{(\mu,\phi)}$ in $(\B_k)^{\times}$. 
Therefore we have 
$$
\align
\Ch[W^{\vep}]\left(\gam^{\mu}\otimes E\right)
&=2^{-r}2^{k-l}(\sqrt{2})^{l-k}
\Ch[W]\left(\sig^{(\mu,\phi)}\otimes E\right)\\
&=(\sqrt{2})^{-2r+k-l)}2^{l}p_{\mu}(x_1,\dots,x_n)\\
&=\left\{
\matrix
\format\l&\quad\l\\
(\sqrt{2})^{l}p_{\mu}(x_1,\dots,x_n)&\t{ if $k$ is even,}\\
&\\
(\sqrt{2})^{1+l}p_{\mu}(x_1,\dots,x_n)&\t{ if $k$ is odd.}\\
\endmatrix\right.
\endalign
$$
On the other hand, using \thetag{1.10}, 
if $k$ is even, then we have 
$$
\align
&\Ch[\bigoplus_{\nu\in DP_k^{+}}
V_{\nu}\sot U_{\nu}\oplus
\bigoplus_{\nu\in DP_k^{-}}
(V_{\nu}\sot U_{\nu})^{+}]
\left(\gam^{\mu}\otimes E\right)\\
&\quad=\sum_{\nu\in DP_k}
\Ch[V_{\nu}](\gam^{\mu})
(\sqrt{2})^{-d(\nu)-l(\nu)}Q_{\nu}(x_1,\dots,x_n)\\
&\quad=\sum_{\nu\in DP_k}
\Ch[V_{\nu}](\gam^{\mu})
(\sqrt{2})^{-\vep(\nu)-l(\nu)}Q_{\nu}(x_1,\dots,x_n)\\
\endalign
$$
and if $k$ is odd, then we have 
$$
\align
&\Ch[\bigoplus_{\nu\in DP_k^{+}}
V_{\nu}\sot U_{\nu}\oplus
\bigoplus_{\nu\in DP_k^{-}}
V_{\nu}\sot U_{\nu}]
\left(\gam^{\mu}\otimes E\right)\\
&\quad=
\sum_{\nu\in DP_k^{+}}
\Ch[V_{\nu}](\gam^{\mu})(\sqrt{2})^{1-l(\nu)}Q_{\nu}(x_1,\dots,x_n)\\
&\quad\hskip2cm+\sum_{\nu\in DP_k^{-}}
\Ch[V_{\nu}](\gam^{\mu})(\sqrt{2})^{-l(\nu)}Q_{\nu}(x_1,\dots,x_n)\\
&\quad=\sqrt{2}\times\sum_{\nu\in DP_k}
\Ch[V_{\nu}](\gam^{\mu})
(\sqrt{2})^{-\vep(\nu)-l(\nu)}Q_{\nu}(x_1,\dots,x_n).\\
\endalign
$$
Consequently the result follows.
\qed\enddemo

\Refs

\pref{1}
{J\. W\. Davies, A\. O\.  Morris }
{The Schur multiplier of the generalized symmetric group}
{J\. London Math\. Soc\. Ser\. 2}{8}{615-620}{1974}

\pref
{2}
{F\. G\. Frobenius}
{\"Uber die Charaktere der symmetrischen Gruppe} 
{  Sitzungsberichte der K\"oniglich Preussischen Akademie der 
Wissenschaften zu Berlin}{}{516-534}{1900}

\pref{3}
{T\. J\'ozefiak}
{Semisimple superalgebras}
{in: Some Current Trends in Algebra, 
Proceedings of the Varna Conference 1986, 
Lecture Notes in Math\.  1352, Springer Berlin}{}{96-113}{1988}

\pref{4}
{T\. J\'ozefiak}
{Characters of projective representations of symmetric groups}
{Expo\. Math\. }{7}{193-247}{1989}

\pref{5}
{T\. J\'ozefiak}
{Schur $Q$-functions and applications}
{Proceedings of the Hyderabad Conference on Algebraic groups, 
Manoj Prakashan}
{}{205-224}{1989}

\pref{6}
{T\. J\'ozefiak}
{A class of projective representations of hyperoctahedral groups and 
Schur $Q$-functions}
{Topics in Algebra, Banach Center Publications, Vol\. 2, Part 2,
PWN-Polish Scientific Publications, Warsaw}
{}{}{1990}

\pref{7}
{V\. G\. Kac}
{Lie superalgebras}
{Adv\. in Math\.}{26}{8-96}{1977}

\pref{8}
{I\. G\. Macdonald}
{Symmetric functions and Hall polynomials}
{Clarendon Press, Oxford}{}{}{1979}

\pref{9}
{B\. E\. Sagan}
{Shifted Tablaeux, Schur $Q$-functions, and a Conjecture of R\. Stanley 
\hskip3cm}
{J\. Combin\. Theory Ser\. A}
{45}{62-103}{1987}

\pref{10}
{I\. Schur}
{\"Uber die Darrstellung der symmetrischen und der alternierenhen 
Gruppe durch gebrohene lineare Substitutionen}{J\. Reine Angew\.  Math\. }
{139}{155-250}{1911}

\pref{11}
{A\. N\. Sergeev}
{Tensor algebra of the identity representation as a module over Lie 
superalgebras $GL(n,m)$ and $Q(n)$}
{Math\.  USSR Sbornik}
{51, No\. 2}
{419-425}{1985}

\pref{12}
{J\. R\. Stembridge}
{Shifted Tableaux and the Projective Representations of Symmetric groups}
{Adv\.  in Math\. }
{74}{87-134}{1989}

\pref{13}
{J\. R\. Stembridge}
{The Projective Representations of the Hyperoctahedral Group}
{J\.  Algebra}
{145}{396-453}{1992}

\pref{14}
{J\. R\. Stembridge}
{Some permutation representations of Weyl groups associated 
with the cohomology of toric varieties}
{Adv\. in Math\. }
{106}{244-301}{1994}

\pref{15}{C\. T\. C\. Wall}{Graded Brauer groups}
{J\. Reine Angew Math\. }{213}{187-199}{1964}

\endRefs

\enddocument